\documentclass{amsart}
\usepackage{amsfonts,amscd}

\newtheorem{theorem}{Theorem}[section]
\newtheorem{lemma}[theorem]{Lemma}
\newtheorem{corollary}[theorem]{Corollary}
\newtheorem{proposition}[theorem]{Proposition}
\theoremstyle{remark}

\theoremstyle{definition}

\numberwithin{equation}{section}
\makeatother

\DeclareMathOperator{\Cdb}{{\mathbb C}}
\DeclareMathOperator{\Rdb}{{\mathbb R}}
\DeclareMathOperator{\Ddb}{{\mathbb D}}
\DeclareMathOperator{\Tdb}{{\mathbb T}}
\DeclareMathOperator{\Ndb}{{\mathbb N}}

\begin{document}

\title[Von Neumann algebraic $H^p$ theory]{Von Neumann algebraic  $H^p$ 
theory}

\date{\today}

\author{David P. Blecher}
\address{Department of Mathematics, University of Houston, Houston, TX
77204-3008}
\email[David P. Blecher]{dblecher@math.uh.edu}
 \author{Louis E. Labuschagne}
\address{Department of Mathematical Sciences, P.O. Box
 392, 0003 UNISA, South Africa}
\email{labusle@unisa.ac.za}
\thanks{To appear in: Proceedings of the Fifth Conference on
Function Spaces", Contemp. Math} 
%Presented at the 5th conference on Function Spaces, and also at
%the Bedlewo conference on noncommutative harmonic
%analysis and probability, 2006.}
\thanks{*Blecher was partially supported by grant DMS 0400731 from
the National Science Foundation.  Labuschagne was partially
supported by a National Research Foundation Focus Area Grant.}

\begin{abstract}
Around 1967, Arveson invented a striking noncommutative
generalization of classical $H^\infty$, known as {\em subdiagonal
algebras}, which include a wide array of examples of interest to
operator theorists. Their theory extends that of the generalized
$H^p$ spaces for function algebras from the 1960s, in an extremely
remarkable, complete, and literal fashion, but for reasons that are
`von Neumann algebraic'.  
%Although there are some  
%of new results, some from a
%dd
%companion work, 
Most of the present paper consists of a
survey of our work on Arveson's algebras, and the attendant $H^p$
theory, explaining some of the main ideas in their proofs, and
including some improvements and short-cuts.   The newest results
%dd
utilize new variants of the noncommutative Szeg\"{o} theorem for
$L^p(M)$, to generalize many of the classical results concerning
outer functions, to the noncommutative $H^p$ context. In doing so we
solve several of the old open problems in the subject.
%dd
 We include full
proofs, for the most part, of the simpler `antisymmetric algebra'
special case of our results on outers.
\end{abstract}

\maketitle

\section{Introduction}

%EE
What does one get if one combines the theory of von Neumann algebras, 
and that of Hardy's $H^p$ spaces? In the 1960s, Arveson suggested
one way in which this could be done \cite{PTAG,AIOA}, via his
introduction of the notion of a
{\em subdiagonal} subalgebra $A$ of a von Neumann algebra $M$.
In the case that $M$ has a finite trace (defined below), $H^p$ may be 
defined to be the closure of $A$ in the noncommutative $L^p$ space
$L^p(M)$, and our article concerns the 
`generalization' of Hardy space results to this 
setting.  In the case that $A = M$, the $H^p$ theory 
collapses to noncommutative
$L^p$ space theory.  At the other extreme, if 
%EE
$A$ contains no selfadjoint elements except scalar multiples of the identity, 
the $H^p$ theory will in the setting where $M$ is commutative, collapse to the 
classical theory of $H^p$-spaces associated to the so-called `weak*
Dirichlet algebras'---a class of abstract function algebras.  
Thus Arveson's setting formally merges noncommutative $L^p$ 
spaces, and the classical theory of $H^p$-spaces for
abstract function algebras.        
%BB
We say more about the latter:
around the early 1960's, it became apparent that many
famous theorems about the classical $H^\infty$ space of
bounded analytic functions on the disk,
could be generalized to the setting of abstract function algebras.
This was the work of 
%B
several notable researchers, and in particular
Helson and Lowdenslager \cite{HL}, and Hoffman \cite{Ho}. The paper
\cite{SW} of Srinivasan and Wang, from the
middle of the 1960s decade, organized and summarized much of this
`commutative generalized $H^p$-theory'.
%BB for abstract function algebras.   
In the last few years, we have shown that all of the results
from this classical survey,
%dd
and essentially everything relevant
in Hoffman's landmark paper \cite{Ho}, extend in a particularly 
literal way to
the noncommutative setting of Arveson's  subdiagonal
subalgebras
% of von Neumann algebras \cite{AIOA}
%EE
(en route solving what seemed to be the major open problems in the
latter subject).   Indeed, as an example of what some might call
`mathematical quantization', or noncommutative (operator algebraic)
generalization of a classical theory, the 
%B
program succeeds to a
degree of `faithfulness to the original'  
which seems to be quite rare. This in turn suggests that the analytic
principles encapsulated in the classical theory 
%B may well be 
are far more
algebraic in nature than
was even anticipated in the 1960s.

%B
Since the classical theories mentioned above are 
%B
quite  beautiful, and since the noncommutative variant is such
a natural place of application of von Neumann algebra
and noncommutative $L^p$ space theory,
this has been a very pleasurable labor, which we are
grateful to have participated in.
%D
In any case, it seems timely to survey the main parts of our work.
We aim this at a general audience, and include a description of a
few of the von Neumann algebraic and noncommutative $L^p$ space
methods that are needed. We also include a somewhat improved route
through some of our proofs. In addition, there are two sections
%dd
describing very recent results from  \cite{BL6}. For example, we
give new variants of the noncommutative Szeg\"{o}'s theorem for
$L^p(M)$, and using these, we generalize many of the classical
results concerning outer functions to the noncommutative $H^p$
context.  In particular, we develop in Section 7 the theory of outer
functions in the
%dd
simple and tidy `antisymmetric algebra' case, proving almost all of
our assertions.  Outers for more general algebras are treated in
\cite{BL6}.

We write $H^\infty(\Ddb)$  for the algebra
of  analytic and bounded functions
on the open unit disc $\Ddb$ in the complex plane.
As we just mentioned, and now describe briefly, around 1960 many
notable mathematicians
attempted to generalize
the theory of $H^p$ spaces.
%BB on the disk to the framework of general function algebras. 
 The setting for their
`generalized $H^p$ function theory' was the following:
Let $X$ be a probability space, and let $A$ be
 a closed unital-subalgebra of $L^\infty(X)$, such that:
\begin{equation} \label{eq1}
\int \, f g = \int \, f \; \int \, g , \qquad f, g \in A .
\end{equation}
The latter was the crucial condition which they had isolated as
underpinning the classical generalized function theory.    Note that
 this clearly holds in the
case of $H^\infty(\Ddb) \subset L^\infty(\Tdb)$, the integral being
the normalised Lebesgue  integral for the unit circle $\Tdb$. We
will suppose that $A$ is weak* closed (otherwise it may be replaced
by its weak* closure). Write $[{\mathcal S}]_p$ for the closure of a
set ${\mathcal S} \subset L^p$ in the $p$-norm, and define $H^p =
[A]_p$ whenever $A$ satisfies any/all of the conditions described in
the following theorem.   Let $A_0 = \{ f \in A : \int \, f = 0 \}$.
Combining fundamental ideas of many researchers, one may then prove
\cite{Ho,SW,HR} that:

\begin{theorem} \label{SW}
For such $A$, the following eight conditions are equivalent:
 \begin{itemize}
\item [(i)]  The weak* closure of
$A + \bar{A}$ is all of $L^\infty(X)$.
\item [(ii)]  The `unique normal state extension property' holds,
that is: if  $g \in L^1(X)$ is nonnegative
with $\int f g =  \int f$
for all $f \in A$, then $g = 1$ a.e..
\item [(iii)] $A$ has `factorization', that is: if $b \in L^\infty(X)$, 
%ll b > 0 is a bit ambiguous
then $b \geq 0$ and is bounded away from 0 iff
$b = |a|^2$ for an invertible $a \in A$.
\item [(iv)]  $A$ is `logmodular', that is:  if $b \in L^\infty(X)$ with 
%ll
$b$ bounded away from 0 and $b \geq 0$,  then
$b$ is a uniform limit of terms of the form $|a|^2$ for an invertible $a \in A$).
\item [(v)]  $A$ satisfies the $L^2$-distance formula in
Szeg\"o's theorem, that is: $\exp \int \log g =  \inf \{ \int
|1-f|^2 g : f \in A , \int f = 0 \}$, for any nonnegative
$g \in
L^1(X)$.
\item [(vi)]  Beurling invariant subspace property: every $A$-invariant
closed subspace of $L^2(X)$ such that $[A_0 K]_2 \neq K$ is
of the form $u H^2$ for a unimodular function $u$.
\item [(vii)]  Beurling-Nevanlinna factorization property, that is: 
every $f \in L^2(X)$ with
%dd
 $\int \, \log |f| > - \infty$ has an `inner-outer
factorization' $f = u h$, with $u$ unimodular and $h \in H^2$ outer
(that is, such that $1 \in [h A]_2$).
\item [(viii)]  Gleason-Whitney property: there is a unique
Hahn-Banach extension to $L^\infty(X)$ of any weak* continuous
functional on $A$, and this extension is weak* continuous.
\end{itemize}
\end{theorem}

It is worth remarking that almost none of the implications here are
clear; indeed the theorem above constitutes a
resum\'{e} of a network connecting
several topics of great interest. The objects characterized here,
the {\em weak* Dirichlet algebras}, are the topic of \cite{SW}.  The
theory goes on to show that they satisfy many other properties that
generalize those of bounded analytic functions on the disk; e.g.\
$H^p$ variants of
properties (vi) and (vii) for all $1 \leq p \leq
\infty$,  Jensen's inequality, the F \& M Riesz theorem, Riesz
factorization, the characterization of outer functions in terms of
$\int \, \log f$,
%dd
the criteria for inner-outer factorization, and so on.

For us, the primary question is
the extent to which all of this survives the passage to
noncommutativity. The remarkable answer is that in the setting of
Arveson's subdiagonal subalgebras of von Neumann algebras,
essentially everything does.   A key point when trying to generalize
this theory to a von Neumann algebra framework, is that one must
avoid most classical arguments which involve exponentials of
functions---since the exponential map behaves badly if the exponent
is not a normal operator. Thus we
avoided some of the later, more
sophisticated, routes through the classical theory (see e.g.\
\cite{Gam}), and went back to the older
more algebraic methods of
Helson-Lowdenslager, Hoffman, and others.  Being based primarily on
Hilbert space methods, these more easily go noncommutative. Although
the {\em statements} of the results which we obtain are essentially
the same as in the commutative case, and although the proofs and
techniques in the noncommutative case may often be modeled loosely
on the `commutative' arguments of the last-mentioned authors, they
usually are much more sophisticated, requiring substantial input
from the theory of von Neumann algebras and noncommutative
$L^p$-spaces.  Sometimes completely new proofs have had to be
invented (as was the case with for example Jensen's formula).

We now review some of the definitions we shall use throughout. For a
set ${\mathcal S}$, we write ${\mathcal S}_+$ for the set $\{ x \in
{\mathcal S} : x \geq 0 \}$.  The word `normal'
applied to linear mappings as usual means
`weak* continuous'. We assume throughout that $M$ is a von Neumann
algebra possessing a faithful normal
tracial state $\tau$.  Here `faithful' means that Ker$(\tau) \cap
M_+ = (0)$, and `tracial' means that $\tau(xy) = \tau(yx)$ for all
$x,y \in M$.  The existence of such $\tau$ implies that $M$ is a
so-called {\em finite} von Neumann algebra. One consequence of this,
which we shall use a lot,
 is that if $x^* x = 1$ in $M$, then $x x^* = 1$ too.
 Indeed $0 = 1 - \tau(x^* x) = \tau(1- x x^*)$, and so $1- x x^* = 0$
 because
$\tau$ is faithful. Applying the above to the partial isometry in
the polar decomposition of any $x \in M$, implies, in operator
theoretic terms, that $x$ is onto iff $x$
is invertible iff $x$ is bounded below.
 From this in turn it follows that for
any $a, b \in M$, $ab$ will be invertible precisely when $a$ and $b$
are separately invertible.
% A similar argument
%shows that $x \in M$ has dense range .

 A {\em tracial subalgebra} of $M$ is a weak* closed subalgebra $A$ of
$M$ such that the 
%B
(unique) trace preserving\footnote{This means that $\tau
\circ \Phi = \tau$.} conditional expectation $\Phi : M \to A \cap
A^* \overset{def}{=} {\mathcal D}$ (guaranteed by \cite[p.\
332]{Tak}) satisfies:
\begin{equation} \label{Eq2}
 \Phi(a_1 a_2) = \Phi(a_1) \,  \Phi(a_2) , \; \; \; a_1, a_2 \in A .
\end{equation}
(Note that (\ref{Eq2}) is a variant of the crucial formula
(\ref{eq1}) underpinning the entire theory.) A {\em finite maximal
subdiagonal algebra} is a tracial subalgebra of $M$ with $A + A^*$
weak* dense in $M$.  For brevity we will usually drop the word
`finite maximal' below, and simply say `subdiagonal algebra'. In the
classical function algebra setting \cite{SW}, one assumes that
${\mathcal D} = A \cap A^*$ is one dimensional, which forces $\Phi =
\tau(\cdot) 1$.  If in our setting this is the case, then we say
that $A$ is {\em antisymmetric}.   
%B
It is worth remarking that the antisymmetric maximal subdiagonal 
subalgebras of commutative von Neumann algebras
are precisely the (weak* closed) weak* Dirichlet algebras.
The simplest example of a noncommutative maximal
subdiagonal algebra is the upper triangular matrices $A$ in $M_n$.
Here  $\Phi$ is the expectation onto the main diagonal.  There are
much more interesting examples from free group von Neumann algebras,
etc.  See e.g.\ \cite{AIOA,Zs,LM,KT,MMS}; and in the
next paragraph we will mention a couple of examples in a
little more detail.  In fact much of Arveson's
extraordinary original paper develops a core of substantial examples
of interest to operator theorists and operator algebraists;
%BB
indeed his examples showed that his theory unified part of 
the existing theory of nonselfadjoint operator algebras.  Note
too that the dropping of the `antisymmetric' condition above, gives the
class of subdiagonal algebras a generality and scope much wider than
that of  weak* Dirichlet algebras.  Thus, for example, $M$ itself is a
maximal subdiagonal algebra (take $\Phi = Id$).  It is also
remarkable, therefore, that so much of the classical $H^p$ theory
does extend to {\em all} maximal subdiagonal algebras.  However the
reader should not be surprised to find some results here which do
require restrictions on the size of ${\mathcal D}$.
%dd For example,
%in Section 7 (resp.\ Section 5) we usually have
%restrictions on the dimension of ${\mathcal D}$.
Truthfully though, in some
of these results the restrictions may well ultimately be able to be
weakened further.

%BB
To get a feeling for how subdiagonal subalgebras can arise, we take a paragraph 
to mention very briefly just two interesting examples.  See the papers referred
to in the last paragraph for more examples, or more details.
 The first of these examples is due to Arveson \cite[Section 3.2]{AIOA}.  Let $G$ be a 
countable discrete group with a
linear ordering  which is invariant under 
left multiplication, say.  For example, any free group is known to have such 
an ordering (see e.g.\ \cite{Neu}).    This implies that 
$G = G_+ \cup G_{-}, G_+ \cap G_{-} = \{1 \}$.  The 
subalgebra generated by $G_+$ in the group von
Neumann algebra of $G$, immediately gives a 
subdiagonal algebra.     For a second example (see \cite{Zs,LM,KT}), if 
$\alpha$ is any one-parameter group of
$*$-automorphisms of a von Neumann algebra $M$ satisfying a
certain ergodicity
property (and in particular, all those arising in the
Tomita-Takesaki theory), naturally gives rise to a subdiagonal
algebra $A \subset M$, coming from those elements of $M$ whose
`spectrum with respect to $\alpha$' lies in the nonnegative 
part of the real line.  

 By analogy with the classical case, we set $A_0 = A \cap
{\rm Ker}(\Phi)$.  For example, if   $A = H^\infty(\Ddb)$ then
 $A_0 = \{ f \in H^\infty(\Ddb) : f(0) = 0 \}$.
For subdiagonal algebras the analogue of $H^p$ is $[A]_p$, the
closure of $A$ in the noncommutative $L^p$ space $L^p(M)$, for $p
\geq 1$. The latter object may be defined to be the completion of
$M$ in the norm $\tau(\vert \cdot \vert^p)^{\frac{1}{p}}$. The
spaces $L^p(M)$ are Banach spaces satisfying the usual duality
relations and H\" older inequalities \cite{FK,PX}.   There is a
useful alternative definition.  For our (finite) von Neumann algebra
$M$ on a Hilbert space $H$, define $\widetilde{M}$ to be the set of
unbounded, but closed and densely defined, operators on $H$ which
are affiliated to $M$ (that is, $T u = u T$ for all unitaries $u \in
M'$). This is a $*$-algebra with respect to the `strong' sum and
product (see Theorem 28 and the example following it in
\cite{Terp}). The trace $\tau$ extends naturally to the positive
operators in $\widetilde{M}$.  If $1 \leq p < \infty$,  then
$L^p(M,\tau) =
 \{a \in \widetilde{M} :  \tau(|a|^p) < \infty\}$,
equipped with the norm $\|\cdot\|_p = \tau(|\cdot|^p)^{1/p}$ (see
e.g.\ \cite{Nel,FK,Terp,PX}). For brevity, we will in the following
write $L^p$ or $L^p(M)$ for $L^p(M,\tau)$.  Note that $L^1(M)$ is
canonically isometrically isomorphic, via an $M$-module map, to the
predual of $M$.   Of course this isomorphism takes $T \in L^1(M)$ to
the normal functional $\tau(T \, \cdot \,)$ on $M$. This isomorphism
also respects the natural positive cones on these spaces (the
natural cone of the predual of $M$ is the space of positive normal
functionals on $M$).

Arveson realized that the appropriate Szeg\"o theorem/formula for
his algebras should read:
$$\Delta(h) = \inf \{ \tau(h |a+d|^2)
: a \in A_0 , d \in {\mathcal D} , \Delta(d) \geq 1 \}$$ for all $h
\in L^1(M)_+$.  Here $\Delta$ is the {\em Fuglede-Kadison
determinant}, originally defined on $M$ by  $\Delta(a) = \exp
\tau(\log |a|)$ if $|a|$ is strictly positive, and otherwise, $\Delta(a) = \inf \,
\Delta(|a| + \epsilon 1)$, the infimum taken over all scalars
$\epsilon > 0$. Classically this quantity of course represents the
geometric mean of elements of $L^\infty$.
%DD
In \cite{BL2} we noted that this  definition  of $\Delta(h)$
 makes  sense
for $h \in L^1(M)$, and this form was used extensively
in that paper.  Although in
\cite{AIOA}, Arveson does define $\Delta$ for normal functionals
(equivalently elements of $L^1(M)$),  the above is not his
original definition, and some work is necessary to prove that the
two are equivalent \cite[Section 2]{BL2}.

%D
In passing,
% although only defined for elements of $L^1(M)$ in
%\cite{BL2},
we remark that the definition above of $\Delta(h)$ makes perfect
sense for $h$ in any $L^q(M)$ where $q > 0$.  Since we will need
this later we quickly
%dd give the details
explain this point in our setting (see also \cite{Brme,HS}), adapting
the argument in the
third paragraph of \cite[Section 2]{BL3}.  We make use of the Borel
functional calculus for unbounded operators applied to the
inequality
$$0 \leq \log t \leq \frac{1}{q}t^q \quad t \in [1,\infty).$$
Notice that for any $0 < \epsilon < 1$, the function $\log t$ is
bounded on $[\epsilon,1]$. So given $h \in L^1(M)_+$ with $h \geq
\epsilon$, it follows that $(\log h) e_{[0,1]}$ is similarly
bounded. Moreover the previous centered equation ensures that $0
\leq (\log h) e_{[1,\infty)} \leq \frac{1}{q}h^q e_{[1,\infty)} \leq
\frac{1}{q}h^q.$ Here $e_{[0,\lambda]}$ denotes the spectral
resolution of $h$. Thus if $h \in L^q(M)$ and $h \geq \epsilon$ then
$\log h \in L^1(M)$.

Unfortunately, the conjectured noncommutative
Szeg\"o formula stated above, and the (no doubt  more important)
associated
{\em Jensen's inequality}
$$\Delta(\Phi(a)) \leq \Delta(a) , \qquad a \in A ,$$
and {\em Jensen formula}
$$\Delta(\Phi(a)) = \Delta(a) , \qquad {\rm invertible} \; a \in A ,$$
resisted proof for nearly 40 years (although Arveson did prove
these for most of the examples that he was interested in).  In 2004,
via a judicious use of
a noncommutative variant of a classical limit formula for the
geometric mean, and a careful choice of recursively defined
approximants,
%DDD
the second author proved in \cite{LL3} that all maximal subdiagonal
algebras satisfy Jensen's formula (and hence the Szeg\"o formula and
Jensen's inequality too by Arveson's work).  Settling this old open
problem opened up the theory to the recent developments surveyed
here.  Of course much of the classical theory had already been
generalized to subdiagonal algebras in Arveson's original and
seminal paper \cite{AIOA}, and in the intervening decades following
it. We mention for example the work of Zsid\'o, Exel, McAsey, Muhly, Saito
(and his school in Japan), Marsalli and West, Nakazi and Watatani,
Pisier, Xu, Randrianantoanina, and others (see our reference list
below, and references therein).  This work, together with the
results mentioned below, yields a complete noncommutative
generalization of all of the classical theory surveyed in \cite{SW}.
Since much of this work has been surveyed recently in \cite[Section
8]{PX}, we will not attempt to survey this literature here.

As a first set of results,
which may be regarded in some sense as a `mnemonic' for
the subject, we obtain the same cycle of theorems as in the classical case:
%, the following result may be viewed as the
%unifying result in the theory.

\begin{theorem} \label{B-L}  For a tracial subalgebra $A$ of
$M$, the following eight conditions are equivalent:
\begin{itemize} \item [(i)] $A$ is maximal subdiagonal, that is:
 $\overline{A + A^*}^{w*} = M$.
  \item [(ii)]  a) \ $L^2$-density of $A + A^*$ in $L^2(M)$; and
b) the unique normal state extension property, that is: if $g \in
L^1(M)_+$,  $\tau(f g) =  \tau(f)$ for all $f \in A$, then $g = 1$.
\item [(iii)] $A$ has factorization, that is:
 an element $b \in M_+$ is invertible iff $b = a^* a$ for an invertible $a
\in A$.
\item [(iv)]  $A$ is logmodular, that is:
if $b \in M_+$ is invertible then $b$ is a uniform limit of terms of
the form $a^* a$ for invertible $a \in A$.
\item [(v)]  $A$ satisfies the Szeg\"o formula above.
\item [(vi)]   Beurling-like invariant subspace condition (described in
Section 4).
\item [(vii)]    Beurling-Nevanlinna factorization property, that is:
%dd (described inSection 4).
%ll (lowercase e)
every $f \in L^2(M)$ such that $\Delta(f) > 0$ has an `inner-outer
factorization' $f = u h$, with $u$ unitary and $h \in H^2$ outer (that is,
$1 \in [h A]_2)$.
\item [(viii)]   Gleason-Whitney property: there is at most one 
normal Hahn-Banach extension to $M$ of any normal functional on $A$.
\end{itemize}
\end{theorem}

It will be noted that there is an extra condition in (ii) that does
not appear in the classical case.  It is interesting that this extra
condition took some years to remove in the classical case (compare
\cite{SW} and \cite{HR}).  Although we have not succeeded yet in
removing it altogether in our case, we have made partial progress in
this direction in  \cite[Section 2]{BL4}.

We have also been able to prove many other generalizations of the
classical generalized function theory in addition to those already
mentioned; for example the F \& M Riesz theorem, $L^p$ versions of
the Szeg\"o formula, the Verblunsky/Kolmogorov-Krein
extension of the Szeg\"{o} formula, inner-outer factorization, etc.
%This is in turn used
In Section 7 we will generalize important aspects of the classical
theory of outer functions to subdiagonal algebras, formally
completing the generalization of \cite{SW}. In this regard we note
that $h \in L^p(M)$ is {\em outer} if $[h A]_p = H^p$. That is, $h
\in H^p$, and $1 \in [h A]_p$.
%dd
This definition is in line with e.g.\ Helson's definition 
of outers in the matrix valued case he considers in 
\cite{Hel} (we thank Q. Xu for this observation). 
 In Section 7 we will restrict our attention to the special case
of antisymmetric subdiagonal algebras, where the theory of
outer functions works out particularly transparently and tidily; the
general case will be
treated in the forthcoming work
\cite{BL6}.   It is worthwhile pointing out however,
 that there are
many interesting antisymmetric maximal subdiagonal algebras besides
the weak* Dirichlet algebras---see \cite{AIOA}.
%dd
Our main results here state 1)\ that $h \in H^p$ is outer iff $\Delta(h)
= \Delta(\Phi(h)) > 0$ (one direction of this is not
quite true in the general case discussed in 
\cite{BL6}); and 2)\ 
%er $A$ is antisymmetric or not) 
if $f \in L^p(M)$ with $\Delta(f) > 0$,
then $f = u h$ for a unitary $u \in M$ and an outer $h \in H^p$
(we prove elsewhere that this is true also in the general, i.e.\ 
non-antisymmetric, case). In
particular, $|f| = |h|$ for this outer $h \in H^p$, which solves an
approximately thirty year old problem (see e.g.\ the discussion in
\cite[p.\ 386]{MMS}, or \cite[Chapter 8]{PX}, particularly lines 8-12
on p.\ 1497 of the latter reference).   We remark that the commutative
case of these results was settled in \cite{Nak}.

We end this introduction by mentioning that there are many other,
more recent, generalizations of $H^\infty$, based around
multivariable analogues of the Sz-Nagy-Foia\c s model theory for
contractions. Many prominent researchers are currently intensively
pursuing these topics, for example Popescu, Arias; Arveson again;
Ball and Vinnikov; Davidson and  Power and their brilliant
collaborators;  and Muhly and Solel. See e.g.\ \cite{Pop,MS} and
references therein.   In essence,  the unilateral shift is replaced
by left creation operators on some variant of Fock space. These
generalizations are very important at the present time, and are
evolving in many directions (with links to wavelets, quantum
physics, conservative linear systems, and so on). Although these
theories also contain variants of parts of the theory of $H^\infty$
of classical domains, 
%B
so that in a superficial reading the endeavours may appear to be similar,
in fact they are quite far removed,
% from the theory of subdiagonal algebras, 
and indeed have nothing in common from 
a practical angle.   For example, those other theories have nothing
to do with (finite) von Neumann algebra techniques, which are absolutely
key for us.  So, for example, if one compares Popescu's theorem
of Szeg\"o type from \cite[Theorem 1.3]{Pop} with the Szeg\"o
theorem for subdiagonal algebras discussed here, one sees that they
are only related in a very formal sense.  
 It is unlikely that the
theory of subdiagonal algebras will merge to any great extent with
these other theories, but certain developments in one theory
might philosophically inspire the other.

\section{Two $L^p$-space tools}

%BB
In this survey we will only be able to prove a selection
of our results, and even then some of the proofs will
be sketchy.
% Many of the proofs we will not be able to explain in great detail;
%BB
This is not the forum for a full blown account, and also some of the
proofs are quite technical.
Nonetheless, it seems worthwhile to explain to a general audience a
couple of the tools, each of which is used several times, and which
are quite helpful in adapting proofs of some classical results
involving integrals, to the noncommutative case.  The first tool is
a useful reduction to the classical case.
This may be viewed as a principle of local commutativity for semifinite
von Neumann algebras which furnishes a link between classical and
noncommutative $L^p$ spaces.
Suppose that $h \in L^1(M)$
is selfadjoint. One may of course view $h$ as a normal functional on
$M$, but we will want instead to view $h$ as an unbounded
selfadjoint operator on the same Hilbert space on which $M$ acts, as we
indicated above.   As we shall show, $h$ may be regarded as a
function in a classical $L^1$ space.
 Let $M_0$ be the von Neumann algebra generated by $h$ (see e.g.\
\cite[p.\ 349]{KR}).  This is a commutative subalgebra of $M$, and
it is the intersection of all von Neumann algebras with which $h$ is
affiliated. Let $\psi = \tau_{\vert M_0}$. Since $\psi$ is a
faithful normal state on $M_0$, it is a simple consequence of the
Riesz representation theorem applied to $\psi$, that $M_0 \cong
L^\infty(\Omega,\mu_\tau)$ $*$-algebraically, for a measure space
$\Omega$ and a Radon probability measure $\mu_\tau$. Also, $L^1(M_0)
\subset L^1(M)$, and $L^1(M_0) \cong L^1(\Omega,\mu_\tau)$.
 Via these identifications, $\tau$ restricts to the
 integral $\int_\Omega \; \cdot \, d
\mu_\tau$ on $L^1(\Omega,\mu_\tau)$, and $h$ becomes a real valued
function in $L^1(\Omega,\mu_\tau)$.
In particular $\Delta(h)$ also survives the
passage to commutativity since it is clear from the above that
$\Delta(h) = \inf_{\epsilon > 0} \int_\Omega \; \log(|h| +
\epsilon 1) \, d \mu_\tau$.
Thus for many purposes, we are now back in
the classical situation.

The second technique we will use is `weighted noncommutative $L^p$
spaces' $L^p(M, h)$.  Here $h \in L^1(M)_+$. We define $L^2(M, h)$
to be the completion of $M$ in the inner product $$\langle a,
b\rangle_h = \tau(h^{1/2}b^*ah^{1/2}), \qquad a,b \in M .$$ Note
that $L^2(M, h)$ can be identified unitarily, and as $M$-modules,
with the closure of $M h^{1/2}$ in $L^2(M)$. Let $a \mapsto \Psi_a$
be the canonical inclusion of $A$ in $L^2(M,h)$. There is a
canonical normal $*$-homomorphism representing $M$ as an algebra of
bounded operators on $L^2(M, h)$.  Indeed define
$$\pi(b)\Psi_a = \Psi_{ba} , \qquad a,b \in M,$$ and then extend
this action to all of $L^2(M, h)$.  This is very closely connected
to the famous notion of the `standard form' or `standard
representation' of a von Neumann algebra (see e.g.\ \cite{Tak2}).

More generally, we define $L^p(M,h)$ to be the completion in
$L^p(M)$ of $M h^{\frac{1}{p}}$.   Note that if $e$ is the support
projection of a positive $x \in L^p(M)$ (that is, the smallest
projection in $M$ such that $e x = x$, or equivalently, the
projection onto the closure of the range of $x$, with $x$ regarded
as an unbounded operator), then it is well known (see e.g.\
\cite[Lemma 2.2]{JS}) that $L^p(M) e$ equals the closure in $L^p(M)$
of $M x$. Hence $L^p(M,h) = L^p(M) e$, where $e$ is the support
projection of $h$. Now for any projection $e \in M$ it is an easy
exercise to prove that the dual of $L^p(M) e$ is
 $e L^q(M)$ (see e.g.\ \cite{JS}).
It follows that  the dual of $L^p(M,h)$ is the variant 
of $L^q(M,h)$ where we consider the completion of 
$h^{\frac{1}{q}} M$.

This procedure corresponds in the classical case, to a Radon-Nikodym
derivative, or to `weighting' a given measure.

\section{The equivalences (i)--(v) in Theorem \ref{B-L}}

In this section we indicate a somewhat simplified route through the
equivalences (i)--(v) in Theorem \ref{B-L} above, which are
originally from \cite{BL2}. For those familiar with \cite{BL2}, we
remark that the approach here 1)\  avoids the use of Lemmas 3.3 and
5.1, and Corollary 4.7, from that paper, 2)\  proves the
implication (ii) $\Rightarrow$ (iii) in the theorem, which hinges on
the property of `$\tau$-maximality' discussed below, more directly,
and 3)\ is more self-contained, avoiding some of the reliance on
results from other papers. We will still need to quote \cite{E} in
one place, and we will need a few facts about the Fuglede-Kadison
determinant from \cite{AIOA}.

Amongst the circle of equivalences in Theorem \ref{B-L}, it is
trivial that (iii) $\Rightarrow$ (iv), and fairly obvious that (i)
$\Rightarrow$ (ii). Indeed if (i) holds, and if $g \in L^1(M)_+$
with $\tau(f g) = \tau(f)$ for all $f \in A$, then $g - 1 \in
A_\perp$. Since $g-1$ is selfadjoint we deduce that $g - 1 \in (A +
A^*)_\perp = (0)$. Similarly, if $g \in L^2(M)$ with $g \perp A +
A^*$, then since $g \in L^1(M)$, we see that $g = 0$.  So (ii) holds.

We describe briefly some of the main ideas in the proof from
\cite{BL2} that (iv) $\Rightarrow$ (v).  This is a slight
generalization of the solution from \cite{LL3} of Arveson's long
outstanding problem as to whether (i) implied (v) and the Jensen
formula/inequality mentioned in Section 1.  Arveson had proved that
(i) implied (iii) (another proof is sketched below,
%dd
which is longer than Arveson's but can be used to yield some other
facts too), and had also noted that if (i) held then (v) was
equivalent to the Jensen inequality or the Jensen's formula. By
means of some technical refinements to these arguments,
\cite[Proposition 3.5]{BL2} proves that the validity of Jensen's
formula, Jensen's inequality, and the Szeg\"o formula are
progressively stronger statements, with all three being equivalent
if $A$ is logmodular. To see that logmodularity indeed does imply
the Szeg\"o formula, one therefore need only adapt the argument from
\cite{LL3} to show that (iv) implies Jensen's formula (cf.
\cite[Proposition 3.1]{BL2}).

 Next we prove that (v) implies the unique normal state
extension property (that is, (ii)b). Suppose that we are given an
$h \in L^1(M)_+$,
 such that $\tau(h a) = \tau(a)$ for all $a \in A$.
Then $\tau(h a) = 0$ for all $a \in A_0$, and hence also for all $a
\in A_0^*$, since $\overline{\tau(h a^*)} = \tau(ha)$. If $a  \in
A_0$, and $d \in {\mathcal D}$, then
$$\tau(h |a+d|^2) =
\tau(h |a|^2 + h d^* a + h a^* d + h |d|^2) = \tau(|a| h |a| +
|d|^2) \geq \tau(|d|^2) .$$ Appealing to the Szeg\"o formula in (v), we
deduce that
$$\Delta(h) \; = \; \inf \{ \tau(|d|^2) : d \in
{\mathcal D}, \Delta(d) \geq 1 \} \; \leq \; 1 .$$ By
\cite[4.3.1]{AIOA}, we have
$$\tau(|d|^2) \geq \Delta(|d|^2) = \Delta(|d|)^2 = \Delta(d)^2 .$$
It follows that
$$\Delta(h) \; = \; \inf \{ \tau(|d|^2) : d \in
{\mathcal D}, \Delta(d) \geq 1 \}   \; = \;  1  .$$ By hypothesis,
we also have $\tau(h) = \tau(1) = 1$. We now reduce to the classical
case as we described earlier in this section, so that the von
Neumann algebra generated by $h$ is $*$-isomorphic to
$L^\infty(\Omega,\mu)$ $*$-algebraically, for a measure space
$\Omega$ and a Radon probability measure $\mu$, and $h$ becomes a
real-valued function in $L^1(\Omega,\mu)$.  We have $\int_\Omega \;
h \, d \mu = 1 = \exp(\int_\Omega \;  \log h \, d \mu)$. It is an
elementary exercise in real analysis to show that this forces $h = 1$
a.e.\  This proves (ii)b.

We now sketch the proof that (iii) implies (i).  This requires three
technical background facts, which we now state.   In
\cite[Proposition 3.2]{BL2} it is shown that (iii) implies that for
$h \in L^1(M)$, we have
$$\Delta(h) \; = \; \inf \{ \tau(|ha|) : a \; \text{invertible in} \, A ,
\Delta(a) \geq 1 \} .$$ This in turn is an extension of facts about
the Fuglede-Kadison determinant from \cite{AIOA}.  Since the proof
is rather long (using the Borel functional calculus for unbounded
selfadjoint operators affiliated to a von Neumann algebra) we omit
the details.
%dd The interested reader will find similar strategies in Corollary
%\ref{vnc2} of this paper.
 Next, one shows that the last displayed formula,
together with Jensen's formula (which we already know to be a
consequence of (iii)), implies that if $h \in L^1(M)$ with $\tau(ha)
= 0$ for every $a \in A$, then $\Delta(1-h) \geq 1$.
See \cite[Lemma 3.3]{BL2} for the short calculation.  The third fact that
 we shall need is that,
just as in the classical case, if $h \in L^1(M)$ is selfadjoint, and
if for some $\delta
> 0$ we have
$$\Delta(1 - th) \geq 1 \; , \qquad t \in (-\delta,\delta),$$
then $h = 0$.  This is a noncommutative version of an extremely
elegant and useful lemma which seems to have been proved for $L^2$
functions by Hoffman \cite[Lemma 6.6]{Ho}, and then extended to
$L^1$ functions by R.  Arens in what apparently was a private
communication to Hoffman.  The noncommutative version is proved by
reducing to the classical case as we described earlier in this
section, so that the von Neumann algebra generated by $h$ is
$*$-isomorphic to $L^\infty(\Omega,\mu)$ $*$-algebraically, for a
measure space $\Omega$ and a Radon probability measure $\mu$, and
$h$ becomes a function in $L^1(\Omega,\mu)$.  We are now back in the
classical situation, and one may invoke the classical result
mentioned above. See \cite[Section 2]{BL2} for more details if
needed.

We now explain how these facts give the implication (iii)
$\Rightarrow$ (i).   To show that $A + A^*$ is weak* dense in $M$,
it suffices
 to show that if $h \in (A + A^*)_\perp$
then $h = 0$.  Since $A + A^*$ is selfadjoint, it is easy to see
that $h \in (A + A^*)_\perp$ if and only if $h + h^* \in (A +
A^*)_\perp$ and $i (h-h^*)
 \in (A + A^*)_\perp$.
We may therefore assume that $h$ is selfadjoint. By the facts in the
last paragraph, $\Delta(1-th) \geq 1$ for every $t \in \Rdb$, and
this implies that  $h = 0$.

We now complete this circle of equivalences by using weighted $L^2$
space arguments (discussed in Section 2) to show that (ii) implies
(iii), and that (v) implies the $L^2$-density of $A + A^*$ in
$L^2(M)$ (that is, (ii)a).

We say that $A$ is {\em $\tau$-maximal} if $A = \{ x \in M : \tau(x
A_0) = 0 \}$.
%dd
The following is essentially due to Arveson (as
pointed out to us by Xu, the proof in \cite{AIOA} of `factorization'
essentially only uses the hypotheses 1 or 2 below).  However we give
a somewhat different proof since the method will be used again
immediately after the theorem.

\begin{theorem} \label{pfac} Let $A$ be a tracial subalgebra of $M$.
Consider the following statements:
\begin{enumerate}
\item $A$ satisfies $L^2$-density and the unique normal state extension
property;
\item $A$ satisfies $\tau$-maximality and the unique normal state
extension property;
\item $A$ has factorization.
\end{enumerate}
The following implications hold: $(1) \Rightarrow (2) \Rightarrow
(3)$.
\end{theorem}

\begin{proof}  (Sketch) \
To see that (1) $\Rightarrow$ (2)
 we will need the space $A_\infty = [A]_2 \cap
M$.
It is an exercise to show  that $A_\infty$ is also a tracial
subalgebra of $M$, with respect to the same expectation $\Phi$ (see
\cite[Theorem 4.4]{BL2}).   However, an adaption of a beautiful von
Neumann algebraic argument of Exel's from \cite{E},
shows that if $A$
satisfies hypothesis (1) then there is no such properly larger
tracial algebra. See the proof of \cite[Theorem 5.2]{BL2}.   Thus $A
= A_\infty$.

By $L^2$-density, $L^2(M) = [A]_2 \oplus [A_0^*]_2$.  Thus if
 $x \in M$ with $\tau(x A_0) = 0$, then $x \in (A_0^*)^\perp
= [A]_2$, so that $x \in A_\infty = A$.  Hence $A$ is
$\tau$-maximal.

The proof that $(2) \Rightarrow (3)$ rests on a  slight modification
of the proof of \cite[Theorem 4.6(a)]{BL2}. If $A$ is  a
$\tau$-maximal tracial subalgebra of $M$, and if $b \in A_\infty$,
then there exists a sequence $(a_n) \subset A$ with $L^2$-limit $b$.
If $c \in A_0$, then $a_n c \to b c$ in $L^2(M)$, and
$$\tau(b c)  \; = \;  \lim_n \tau(a_n c) \; = \; 0 .$$
Thus $b \in A$.  Therefore $A_\infty = A$. Now suppose that in
addition $A$ satisfies the unique normal state extension property.
If  $b \in M_+$ is invertible,  we consider the weighted
noncommutative $L^2$ spaces' $L^2(M, b)$. Let $p$ be the orthogonal
projection of $1$ onto the subspace $[A_0]_2$, taken with respect to
the weighted inner product.  In \cite[Theorem 4.6(a)]{BL2} it is shown
 that $(1 - p) b (1 - p^*) \in
L^1({\mathcal D})_+$, and that also $(1 - p) b (1 - p^*)  \geq
\epsilon 1$, for some $\epsilon > 0$. Thus this element has a
bounded inverse
in $\mathcal{D}$.  Set $e = ((1 - p) b (1 - p^*))^{-\frac{1}{2}} \in
{\mathcal D}$, and let $a = e(1-p) \in [A]_2$.   It is routine to
see that $a$ is bounded, so that $a \in M$. Hence $a \in M \cap
[A]_2 = A_\infty = A$. Since $1 = a b a^*$, and since $M$ is a
finite von Neumann algebra, we also have $1 = b a^* a$, so that
$b^{-1} = |a|^2$.

For any $a_0 \in A_0$ one sees that
$$\tau(a^{-1} a_0) = \tau(b(1-p^*) e a_0)
= 0, $$ since $ea_0 \in [A_0]_2$, and $1-p \perp [A_0]_2$ in the
weighted inner product. Thus
$$a^{-1} \in \{x \in M: \tau(xA_0) = 0\} = A,$$
using $\tau$-maximality.   We deduce that $A$ has factorization.
\end{proof}

Finally, we say a few words about the tricky implication (v) implies
 (ii)a; full details are given in \cite[Proof of Theorem 4.6
(b)]{BL2}. We suppose that $k \in L^2(M)$ is such that $\tau(k (A +
A^*)) = 0$. We need to show that $k = 0$. Since $A + A^*$ is a
self-adjoint subspace of $M$, we may assume that $k = k^*$. Then $1
- k \in L^1(M)$, so that by an equivalent form of (v), given
$\epsilon > 0$ there exists an invertible element $b \in M_+$ with
$\Delta(b) \geq 1$ and $\tau(|(1 - k) b|) < \Delta(1 - k) +
\epsilon$.
In this case it is a bit more complicated, but one can
modify the weighted $L^2$-space argument
in the second half of the proof of Theorem \ref{pfac}, with $b$ replaced by
$b^{-2}$, to one find an element $a \in A_\infty$ with $b^2 = a^* a$.
This element $a$ is used to prove that
$\Delta(1-k) \geq 1$ (we omit the details).  Replacing $k$ by $tk$,
 where $t \in [-1,1]$,
we conclude that $\Delta(1-tk) \geq 1$ for such $t$. Thus by the
Arens-Hoffman lemma (see the discussion surrounding the centered equation
immediately preceding Theorem \ref{pfac}), we have that $k = 0$ as required.
Hence $A + A^*$ is norm-dense in $L^2(M)$.

\medskip

We end this section with a brief remark concerning algebras with the
unique normal state extension property (that is, (ii)b in Theorem
\ref{B-L}), in hope that they (perhaps in conjunction with results
in \cite[Section 2]{BL4}) lead to a resolution of the question as to
whether (ii)a is really necessary in the list of equivalences in
Theorem \ref{B-L}. The hope is to prove that (ii)b $\Rightarrow$
(ii)a. Now by \cite[Theorem 4.6(a)]{BL2} we know that if $A$
satisfies (ii)b of Theorem \ref{B-L}, then $A_\infty$ has `partial
factorization'; that is if $b \in M_+$ is invertible, then we can
write $b = |a| = |c^*|$ for some elements $a, c \in A$ which are
invertible in $M$, with in addition $\Phi(a) \Phi(a^{-1}) = \Phi(c)
\Phi(c^{-1}) = 1$. Then of course $1 = b^2 b^{-2} = a^* a c c^*$, and
so $a c c^* a^* = 1$. That is, $ac$  is unitary.  If one could show
that it is possible to select $a$ and $c$ in such a way that
$\Phi(ac) = 1$, it would then follow that $1 = ac$.  To see this
simply compute $\tau(|1- ac|^2)$. Thus $A_\infty$ would then have
the factorization property (iii), which would ensure the density of
$A_\infty + A_\infty^*$ in $L^2(M)$. Clearly this implies that
 $A + A^*$ is dense too.
In fact even if all we could show is that $a, c$ can be selected so
that $\Phi(ac)$ is unitary, we could then replace $a$ by
$\widetilde{a} = \Phi(ac)a$, and argue as before to see that
$\widetilde{a}c = 1$, hence that $A_\infty$  has the factorization
property, and hence $A + A^*$ is dense in $L^2(M)$.

\section{Beurling's invariant subspace theorem}
%dd and Beurling-Nevanlinna factorization}

We discuss (right) $A$-invariant subspaces, that is
closed subspaces $K \subset L^p(M)$ with $K A \subset K$. For the sake of
brevity we will therefore suppress the word \emph{right} in the following.
%Le
These have been studied by many authors
e.g.\ \cite{MMS,MW,N,PX,Sai,Zs}, with an eye
to Beurling-type invariant subspace theorems, usually in
the case that $p = 2$.   We mention just two
facts from this literature: Saito showed in \cite{Sai} that any
$A$-invariant subspace of
$L^p(M)$ is the closure of the bounded elements which it contains.
In \cite{N}, Nakazi and Watatani
decompose any $A$-invariant subspace
of $L^2(M)$ into three
orthogonal pieces which they called types  I, II, and III.
In the case that the center
of $M$ contains the center of ${\mathcal D}$, they proved
that every type I
invariant subspace of $L^2(M)$ is of the form $u H^2$ for a
partial isometry $u$.   This is a generalization of the
classical Beurling invariant subspace theorem.

Our investigation into $A$-invariant subspaces of $L^p(M)$
revealed the fact that the appropriate Beurling theorem and
`Wold decomposition' for such spaces, are intimately connected
with the {\em $L^p$-modules} developed recently by
Junge and Sherman \cite{JS}, and their  natural `direct sum',
known as the `$L^p$-column sum'.   We explain these terms:
If $X$ is a closed subspace of $L^p(M)$,
and if  $\{ X_i : i \in I \}$ is
a collection of closed subspaces of $X$, which together densely
span $X$, with the
property that $X_i^* X_j = \{ 0 \}$ if $i \neq j$, then
we say that $X$ is the {\em internal column
$L^p$-sum} $\oplus^{col}_i \, X_i$.  There is also
an extrinsic definition of this sum (see e.g.\
the discussion at the start of Section 4 of \cite{BL3}).
 An $L^p({\mathcal D})$-module is a right ${\mathcal D}$-module
with an $L^{\frac{p}{2}}({\mathcal D})$-valued inner product,
satisfying axioms  resembling those of a Hilbert space or Hilbert
$C^*$-module \cite[Definition 3.1]{JS}. 

%B
The key point, and the nicest feature of
$L^p({\mathcal D})$-modules, is that they may all be written as a
column $L^p$-sum of modules with cyclic vectors, each summand of a
very simple form. It is this simple form that gives the desired
Beurling theorem.  We recall that a vector $\xi$ in a right
${\mathcal D}$-module $K$ is {\em cyclic} if $\xi {\mathcal D}$ is
dense in $K$, and {\em separating} if $d \mapsto \xi d$ is
one-to-one.   
%B
Thinking of our noncommutative $H^p$ theory as a  simultaneous
generalization of noncommutative $L^p$-spaces, and of classical 
%EE
$H^p$ spaces, it is in retrospect not
surprising that the Beurling invariant subspace
classification theorem, in
our setting, should be connected to Junge and Sherman's classification 
of $L^p$-modules.

For our purposes, we prefer to initially decompose $A$-invariant
subspaces of $L^p(M)$ into two and not three summands, which we call
type 1 (= type I) and type 2, and which are defined below. This is
intimately connected to the famous {\em Wold decomposition}. For
motivational purposes we discuss the latter briefly. For example,
suppose that $u$ is an isometry on a Hilbert space $H$, and that $K$
is a subspace of $H$ such that $u K \subset K$. If $A$ is the unital
operator algebra generated by $u$, and $A_0$ the nonunital operator
algebra generated by $u$, then $K$ is $A$-invariant.   (The reader
may keep in mind the case where $u$ is multiplication by the
monomial $z$ on the circle $\Tdb$; in this case $A =
H^{\infty}(\Ddb)$, and $A_0$ is the algebra of functions vanishing
at $0$.) The subspace $W = K \ominus u K = K \ominus [A_0 K]$ is
{\em wandering} in the classical sense that $u^n W \perp u^m W$ for
unequal nonnegative integers $n, m$.  The Wold decomposition writes
 $K = K_1 \oplus K_2$, an orthogonal direct sum, where
$K_1 = [A W]_2$, a condition which matches what is called `type 1'
below.  Also, $u K_2 = K_2$, which is equivalent to $[A_0 K_2]_2 = K_2$,
a criterion which matches what is called `type 2' below.
We have  $u$ 
unitary on $K_2$, whereas the restriction of $u$ to subspaces of
$K_1$ is never unitary (see e.g.\ \cite[Lemma 1.5.1]{Nik}).
The match with the definitions in the next paragraph 
 is exact in the case where $A = H^{\infty}(\Ddb)$.

If $K$ is an $A$-invariant subspace of
$L^2(M)$ then the {\em right wandering subspace} of
$K$ is defined to be $W = K \ominus [K A_0]_2$  (see \cite{N}).
As above,  $K$ is {\em type 1} if the right wandering subspace generates $K$
(that is, $[WA]_2 = K$), and {\em type 2} if  the right wandering subspace is trivial
(that is, $W = (0)$).
For the case $p \neq 2$
we similarly say that $K$ is type 2 if $[K A_0]_p = K$, but
we define the right wandering subspace $W$ of $K$, and type 1 subspaces, a
little differently (see \cite{BL3}). We omit the details in the case
$p \neq 2$, but to ease the reader's
mind, we point out that there is a very explicit
type-preserving lattice isomorphism
between the closed (weak*-closed, if $p = \infty$)
right $A$-invariant subspaces of $L^p(M)$ and those of $L^2(M)$ (see
\cite[Lemma 4.2 \& Theorem 4.5]{BL3}).   Thus the theory of
$A$-invariant subspaces of $L^p(M)$ relies on first
achieving a clear understanding of the $p=2$ case.

\begin{theorem} \label{main}  {\rm \cite{BL3} } \
If $A$ is a maximal subdiagonal
subalgebra of $M$, if $1 \leq p \leq \infty$,
and if $K$ is a closed
(indeed weak* closed, if $p = \infty$) right  $A$-invariant subspace of $L^p(M)$,
then: \begin{itemize}
\item [(1)]   $K$ may be written uniquely as
an (internal) $L^p$-column sum $K_1 \oplus^{col} K_{2}$
of a type 1 and a type 2 invariant subspace of $L^p(M)$, respectively.
\item [(2)]  If $K  \neq (0)$ then $K$ is type 1 if and only if
$K = \oplus_i^{col} \, u_i \, H^p$, for $u_i$  partial isometries in
$M$ with mutually orthogonal ranges and $|u_i| \in {\mathcal D}$.
\item [(3)]
The right wandering subspace $W$ of $K$
is an $L^p({\mathcal D})$-module in the sense of Junge and Sherman
(indeed $W^* W \subset L^{p/2}({\mathcal D})$).
\end{itemize}
Conversely, if $A$ is a tracial subalgebra of $M$ such that
every $A$-invariant subspace of $L^2(M)$
satisfies {\rm (1)} and {\rm (2)} (resp.\ {\rm (1)} and {\rm (3)}),
then $A$ is maximal subdiagonal.
  \end{theorem}

We show now why this is a generalization of Beurling's theorem.
Indeed the proof of (a) below proves a classical generalization of
Beurling's theorem in just four lines.  We recall that $h$ is outer
in $H^p$  iff $[h A]_p = H^p$.

\begin{proposition}  \label{BNpa} Let $A$ be an
antisymmetric subdiagonal subalgebra of $M$, and let $1 \leq p < \infty$.
Then
\begin{itemize} \item [{\rm (a)}]
Every right invariant subspace $K$ of $L^2(M)$ such that
$[K A_0]_p \neq K$, is of the
form $u [A]_p$, for a unitary $u$ in $M$,
\item [{\rm (b)}]
Whenever  $f \in L^p(M)$ with $f \notin [f A_0]_p$,
then $f = uh$, for a unitary $u$ and an outer
$h \in H^p$,
\end{itemize}
\end{proposition}

\begin{proof}   (a) \ The subspace $K$ here is not type 2,
and it follows from Theorem \ref{main} that
$K_1 \neq (0)$, and that there
is a nonzero partial isometry $u \in K_1$ with $|u| \in {\mathcal D}
= \Cdb 1$.  Hence $u^* u = 1$, and so
$u$ is a unitary in $M$.  Thus $K_1 = u H^p$.
Since $K_1^* K_{2} = (0)$, we have $K_{2} = (0)$, and so $K = u H^p$.

(a) $\Rightarrow$ (b) \ By (a),
$[f A]_p = u H^p$ for a unitary
$u \in M$. Thus $f = uh$ for $h \in H^p$.  Also, $h$ is outer, since
$u [h A]_p = [f A]_p = u H^p$, so that $[h A]_p = H^p$.
\end{proof}

Admittedly, our formulation of Beurling's invariant subspace theorem
above is more complicated if $A$ is not antisymmetric. We remark
that it can be shown that an $A$-invariant subspace $K$ is of the
form $u H^p$ for a unitary $u \in M$, if and only if the right
wandering subspace of $K$ is a so-called `standard' representation
of ${\mathcal D}$, or equivalently it has a nonzero separating and
cyclic vector, for the right
action of ${\mathcal D}$ \cite[Corollary 1.2]{BL3}. 
We recall that a vector $\xi$ in a right
${\mathcal D}$-module $K$ is {\em cyclic} if $\xi {\mathcal D}$ is
dense in $K$, and {\em separating} if $d \mapsto \xi d$ is
one-to-one.  Or,
%L
alternatively, appealing to Nakazi and Watatani's result mentioned
above, it will follow that if the center of ${\mathcal D}$ is
contained in the center of ${\mathcal M}$, then every type 1
$A$-invariant subspace of $L^p(M)$ is of the form $u H^p$ for a
partial isometry $u \in M$.

%unitary in $M \cap [f A]_p$. \end{theorem}
Part (b) of Proposition  \ref{BNpa}  is called {\em
Beurling-Nevanlinna factorization}.
%dd For antisymmetric algebras one
In fact one can actually give a precise formula for the $u$ and the
$h$
%dd
(see \cite{BL6} and Corollary \ref{pth2}   below). There are also
 generalized Beurling-Nevanlinna factorizations
 in the case that $A$ is not antisymmetric, which we shall discuss
 in Section 7.

\section{The F. and M. Riesz and
Gleason-Whitney theorems}

The classical form of the F.\ and M.\ Riesz theorem 
states that if $\mu \in C(\Tdb)^*$ is a  complex 
measure on the circle whose Fourier coefficients vanish 
on the negative integers (that is, 
%EE
$\mu$ annihilates the trigonometric polynomials
$e^{i n \theta}$, for $n \in \Ndb$), then $\mu$ is
absolutely continuous with respect to Lebesgue measure
(see e.g.\
p.\ 47 of \cite{Ho}).
This is known to fail for weak* Dirichlet algebras; and
hence it will fail for subdiagonal algebras too.  However there is
an equivalent version of the theorem which is true for weak*
Dirichlet algebras \cite{Ho,SW}, and it holds too for a maximal
subdiagonal algebra $A$ in $M$, with ${\mathcal D}$ finite
dimensional.  Moreover, one can show that this dimension condition
is necessary and sufficient for the theorem to hold (see
\cite{BL4}).

\begin{theorem} \label{F&M}  {\rm (Noncommutative
 F. and M. Riesz theorem)} \ If $A$ is as above, and if a
functional $\varphi \in M^*$ annihilates $A_0$ then the normal part
of $\varphi$ annihilates $A_0$, and the singular part annihilates
$A$.
\end{theorem}

Exel proved a striking `norm topology' variant in 1990 \cite{E2},
but unfortunately it does not apply to subdiagonal algebras, which
involve some extra complications.
The strategy of the proof is to translate the main ideas of the
classical proof to the noncommutative context by means of a
careful and lengthy Hilbert space analysis of the `GNS' construction
for a subdiagonal algebra, using the weighted $L^2$ spaces which we
discussed in Section 2, noncommutative Lebesgue-Radon-Nikodym
decomposition of states, etc.

The F. and M. Riesz theorem, like its classical counterpart,
has many applications.  The
main one that we discuss here is the Gleason-Whitney
theorem, and the equivalence with (viii) in Theorem \ref{B-L}.

If $A \subset M$ then we say that  $A$ has property (GW1) if every
extension to $M$ of any normal functional on $A$, keeping the same
norm, is normal on $M$. We say that  $A$ has property (GW2) if there
is at most one normal extension to $M$ of any normal functional on
$A$, keeping the same norm. We say that $A$ has the {\em
Gleason-Whitney property} (GW) if it possesses (GW1) and (GW2). This
is simply saying that there is a unique extension to $M$ keeping the
same norm of any normal functional on $A$, and this extension is
normal.

\begin{theorem}  \label{GW}  If $A$ is a tracial
subalgebra of $M$ then $A$ is maximal subdiagonal if and only if it
possesses property (GW2). If ${\mathcal D}$ is finite dimensional,
then $A$ is maximal subdiagonal if and only if it possesses property
(GW).  \end{theorem}

\begin{proof}  We will simply prove that
if ${\mathcal D}$ is finite dimensional, and if $A$ is maximal
subdiagonal, then $A$ satisfies (GW1).  See \cite[Theorem 4.1]{BL4}
for proofs of the other statements.

Let $\rho$ be a norm-preserving extension of a normal functional
$\omega$
% = \tau(f \cdot)$
on $A$.  By basic functional analysis,
 $\omega$ is the restriction of  a
normal functional $\tilde{\omega}$ on $M$.
% Here $f \in L^1(M)$.
  We may write $\rho = \rho_n + \rho_s$,
where $\rho_n$ and $\rho_s$ are respectively the normal and singular
parts, and $\Vert \rho \Vert = \Vert \rho_n  \Vert + \Vert \rho_s
\Vert$.
%Suppose that $\rho_n = \tau(\tilde{f} \cdot)$ on $M$.
Then $\rho - \tilde{\omega}$
%\tau(f \cdot)$
annihilates $A$, and hence by our F.\ and M.\ Riesz theorem both the
normal and singular parts, $\rho_n - \tilde{\omega}$ and $\rho_s$
respectively, annihilate $A_0$.  Hence they annihilate $A$, and in
particular $\rho_n = \omega$ on $A$.  But this implies that
$$\Vert \rho_n \Vert
+ \Vert \rho_s \Vert = \Vert \rho \Vert = \Vert \omega  \Vert \leq
\Vert \rho_n \Vert .$$
We conclude that $\rho_s  = 0$.   Thus $A$ satisfies (GW1).
\end{proof}

There is another (simpler) variant of the Gleason-Whitney theorem
\cite[p.\ 305]{Ho}, which is quite easy to prove from the F. \& M.
Riesz theorem:

\begin{theorem}  \label{GW2}  Let  $A$ be a maximal subdiagonal subalgebra of $M$
with ${\mathcal D}$ finite dimensional. If $\omega$ is a normal
functional on $M$ then $\omega$ is the unique Hahn-Banach extension
of its restriction to $A + A^*$.  In particular, $\Vert \omega \Vert
= \Vert \omega_{\vert A + A^*} \Vert$ for any such $\omega$.
\end{theorem}

\begin{corollary}  \label{Kap} {\rm
(Kaplansky density theorem for subdiagonal algebras) } \
 Let  $A$ be a maximal subdiagonal subalgebra of $M$
with ${\mathcal D}$ finite dimensional. Then the unit ball of $A +
A^*$ is weak* dense in ${\rm Ball}(M)$.
\end{corollary}

\begin{proof}  If $C$ is the unit ball of $A + A^*$, it follows from the
last theorem  that the pre-polar of $C$ is ${\rm Ball}(M_\star)$. By
the bipolar theorem, $C$ is  weak* dense in ${\rm Ball}(M)$.
\end{proof}

\section{The Fuglede-Kadison determinant and Szeg\"o's theorem for $L^p(M)$}

In this section $A$ is a maximal subdiagonal algebra in $M$.

We  defined the Fuglede-Kadison determinant for elements of $L^q(M)$
in Section 1, for any $q > 0$.  In
%dd
\cite{Brme,HS} it is proved that this determinant has the following
basic properties, which are used often silently in the next few
sections.
%For convenience we include a proof:

\begin{theorem} \label{Fkd}  If $p > 0$ and $h \in L^p(M)$  then
\begin{itemize}
\item [(1)]  $\Delta(h) = \Delta(h^*) = \Delta(|h|)$.
\item [(2)]  If $h \geq g$ in $L^p(M)_+$ then $\Delta(h) \geq \Delta(g)$.
\item [(3)]  If $h \geq 0$ then $\Delta(h^q) = \Delta(h)^q$ for any $q > 0$.
\item [(4)]  $\Delta(h b) = \Delta(h)\Delta(b) = \Delta(b h)$ for any
$b \in L^q(M)$ and any $q > 0$.
\end{itemize}
\end{theorem}

%ADD  proof from BL6, and comment after

%{\bf Remark.}   We take the opportunity here to use a fact from the
%last proof to correct  what appears to be a small gap early in the
%proof of \cite[Theorem 2.1]{BL2}.  This  occurs around at the point
%where it is proved that $\Delta(h)$ is $\leq$ the infimum in (2.1)
%there. One can fill the gap by replacing $h e_{[0,n)}$ in parts of
%that calculation by what is called $\widetilde{h_n}$  above; and use
%the fact, proved above, that $\Delta(\widetilde{h_n}) \to \Delta(h)$.

%\medskip

 Throughout the rest of this section, $A$ is a maximal
subdiagonal algebra in $M$. It is proved in \cite{BL4} that for
 $h \in L^1(M)_+$ and $1 \leq p < \infty$, we have
$$\Delta(h) = \inf \{ \tau(h |a+d|^p) : a \in A_0 , d \in {\mathcal
D} , \Delta(d) \geq 1 \} .$$ A perhaps more useful variant of this
formula is as follows:

\begin{theorem} \label{vnc2} {\rm \cite{BL6} } \  If $h \in L^q(M)_+$
and $0 < p, q < \infty$,
we have $$\Delta(h) = \inf \{\tau(|h^{\frac{q}{p}} a|^p)^{\frac{1}{q}}
: a \in A, \Delta(\Phi(a)) \geq 1 \} = \inf
\{\tau(|a h^{\frac{q}{p}}|^p)^{\frac{1}{q}} : a \in A, \Delta(\Phi(a))
\geq
 1 \}.$$
The infimums are unchanged if we also require $a$ to be invertible
in $A$, or if we require $\Phi(a)$ to be invertible in ${\mathcal
D}$.
\end{theorem}

\begin{corollary} \label{vncc2}  {\rm (Generalized Jensen inequality)
\; \cite{BL6} } \
 Let $A$ be a maximal subdiagonal algebra. For any $h \in H^1$
 we have $\Delta(h) \geq \Delta(\Phi(h))$.
  \end{corollary}

%LL Is this OK? the book of Maddox -- Elements of Funct. Anal. (2nd ed.) --
% discusses the commutative variants
%For the uninitiated we mention
We recall that although $L^p(M)$ is not a normed space
if $1 > p > 0$, it is a so-called linear metric space with metric given by
$\|x - y\|_p^p$ for any $x, y \in L^p$ (see \cite[4.9]{FK}). Thus although
the unit ball may not be convex, continuity still respects all elementary
linear operations.

\begin{corollary} \label{vnc3}  Let $h \in L^q(M)_+$ and $0 < p, q < \infty$.
If $h^{\frac{q}{p}} \in [h^{\frac{q}{p}} A_0]_p$, then
$\Delta(h) = 0$.  Conversely, if $A$ is antisymmetric
and $\Delta(h) = 0$, then $h^{\frac{q}{p}} \in [h^{\frac{q}{p}} A_0]_p$.
Indeed if $A$ is antisymmetric, then
$$\Delta(h) = \inf \{\tau(|h^{\frac{q}{p}} (1 - a_0)|^p)^{\frac{1}{q}} :
a_0 \in A_0 \} .$$
\end{corollary}

\begin{proof}  The first assertion follows by
taking $a$ in the infimum in Theorem \ref{vnc2} to be of the form $1
- a_0$ for $a_0 \in A_0$.

If $A$ is antisymmetric, and if $t \geq 1$ with
$\tau(|h^{\frac{q}{p}} (t1 + a_0)|^p)^{\frac{1}{q}} < \Delta(h) +
\epsilon$, then $\tau(|h^{\frac{q}{p}} (1 + a_0/t)|^p)^{\frac{1}{q}}
< \Delta(h) + \epsilon$. From this the last assertion follows that
the infimums in Theorem \ref{vnc2} can be taken over terms of the
form $1 + a_0$ where $a_0 \in A_0$. If this infimum was $0$ we could
then find a sequence $a_n \in A_0$ with $h^{\frac{q}{p}} (1 + a_n)
\to 0$ with respect to $\Vert \cdot \Vert_p$. Thus $h^{\frac{q}{p}}
\in [h^{\frac{q}{p}}A_0]_p$.
\end{proof}

{\bf Remark.}  The converse in the last result is false for general
maximal subdiagonal algebras (e.g.\ consider $A = M = M_n$).

\medskip

The classical strengthening of Szeg\"o's formula, to the case of
general positive linear functionals, extends even to the
noncommutative context.  Although the classical version of this
theorem is usually attributed to Kolmogorov and Krein, we have been
informed by Barry Simon that Verblunsky proved it in the mid
1930's (see e.g.\ \cite{Ver}):

\begin{theorem} \label{SKK}  {\rm \cite{BL6} \ (Noncommutative
Szeg\"o-Verblunsky-Kolmogorov-Krein theorem)} \ Let
$\omega$ be a positive linear
functional on $M$, and let $\omega_n$ and $\omega_s$
be its normal and singular parts
respectively, with
$\omega_n = \tau(h \, \cdot \, )$ for $h \in L^1(M)_+$.
If $\mathrm{dim}(\mathcal{D}) < \infty$, then
$$\Delta(h) = \inf \{\omega(|a|^2) : a \in A, \Delta(\Phi(a)) \geq 1 \}.$$
%dd this is redundant: The infimum remains unchanged if we also require
%$\Phi(a)$ to be invertible in ${\mathcal D}$.
\end{theorem}

%EE
After seeing this result, Xu was able to use our 
Szeg\"o formula, and facts about singular states, to remove the 
hypothesis that $\mathrm{dim}(\mathcal{D}) < \infty$, and 
to replace the `$2$' by a general $p$. We briefly sketch 
Xu's proof for the case $p = 2$. Firstly note that the fact 
that $A$ has factorization forces the sets $\{|a|^2 : a \in A, 
\Delta(\Phi(a)) \geq 1 \}$ and $\{x : x \in M^{-1} \cap M_+, \Delta(x)) 
\geq 1 \}$ to have a common closure. Hence for any continuous 
linear functional $\rho$ on $M$, we obtain  $$\inf 
\{\rho(|a|^2) : a \in A, \Delta(\Phi(a)) \geq 1 \} = \inf \{\rho(x) 
: x \in M^{-1} \cap M^+, \Delta(x)) \geq 1 \}.$$(See for example 
\cite[4.4.3]{AIOA}.) It is clear that $$\inf \{\omega_n(|a|^2) 
: a \in A, \Delta(\Phi(a)) \geq 1 \} \leq \inf \{\omega(|a|^2) : 
a \in A, \Delta(\Phi(a)) \geq 1 \} .$$ Hence if we can show that 
$\inf \; \omega_n(x)  \geq \inf \; \omega(x)$, these infimums taken over 
the set of $x \in M^{-1} \cap M^+, \Delta(x)) \geq 1$,
then the previously centered equality combined with the Szeg\"o formula, 
will ensure that $\Delta(h) = \inf \{\omega(|a|^2) : a \in A, 
\Delta(\Phi(a)) \geq 1 \}$ as required. By \cite[Theorem III.3.8]{Tak} 
there exists an increasing net $(e_i)$ of projections in the kernel of 
$\omega_s$ such that $e_i \rightarrow 1$ strongly. Given $\epsilon > 0$, 
set $$x_i = \epsilon^{\tau(e_i)-1}(e_i + \epsilon e_i^{\perp}).$$Clearly 
$x_i \in M^{-1} \cap M^+$. Moreover by direct computation $\Delta(x_i) 
= 1$. Now let $x \in M^{-1} \cap M^+$ be given with $\Delta(x) \geq 1$. 
The strong convergence of the $e_i$'s to 1 ensures that $x_ixx_i$ 
converges to $1x1 = x$ in the weak* topology. The weak* continuity of 
$\omega_n$ then further ensures that 
$$\lim\sup\omega(x_ixx_i) = \omega_n(x) + \lim\sup\omega(x_ixx_i) 
\leq \omega_n(x) + \epsilon^2 \|x\|\omega_s(1).$$(Here we have made use 
of the facts that $\omega_s(e_i) = 0$, and that $x_ixx_i$ is
dominated by $\|x\|\epsilon^{2(\tau(e_i)-1)}(e_i + \epsilon^2 e_i^{\perp})$.) Since 
$\epsilon > 0$ is arbitrary, and since $\lim\sup\omega(x_ixx_i)$ dominates 
the infimum of $\omega(x)$, for $x \in M^{-1} \cap M^+, \Delta(x)) \geq 1 \}$, the 
required inequality follows.

\section{Inner-outer factorization and the characterization of
outers}

This section is entirely composed of
%dd
very recent results, and we include almost all the proofs.  In most
of the section $A$ is an antisymmetric maximal subdiagonal algebra;
the much more complicated general case is discussed in more detail
in \cite{BL6}. We recall that if $h \in H^1$ then $h$ is {\em outer}
if $[h A]_1 = H^1$.   An {\em inner} is a unitary which happens to
be in $A$.  We assume $p \geq 1$ throughout this section.

\begin{lemma} \label{L3} Let $1 \leq p \leq \infty$,
and let $A$ be a maximal subdiagonal algebra. Then $h \in L^p(M)$
and
 $h$ is outer in $H^1$, iff
$[h A]_p = H^p$.  (Note that $[\cdot]_\infty$ is the weak* closure.)
%That is, the elements of $L^p(M)$ which are also
%outers in $H^1$ are precisely the $h \in H^p$

If these hold, then $h \notin [hA_0]_p$.  \end{lemma}

\begin{proof} It is
obvious that if $[h A]_p = H^p$ then $[h A]_1 = H^1$. Conversely, if
$[h A]_1 = H^1$ and $h \in L^p(M)$, then by \cite[Lemma 4.2]{BL3}
(the proof of the assertion we are using works for all $p$)
and \cite[Proposition 2]{Sai}, we have
$$[h A]_p = [h A]_1 \cap L^p(M) = H^1 \cap L^p(M) = H^p.$$

If $h \in [hA_0]_p$ then $1 \in [h A]_p \subset [[h A_0]_p A]_p
\subset [hA_0]_p$. Now $\Phi$ continuously extends to a map which
contractively maps $L^p(M)$ onto $L^p({\mathcal D})$ (see e.g.\
Proposition 3.9 of \cite{MW}).  If $h a_n \to 1$ in $L^p$, with $a_n
\in A_0$, then
$$0 =  \Phi(h) \Phi(a_n) = \Phi(h a_n) 
\to \Phi(1) = 1 ,$$
This forces $\Phi(1) = 0$, a contradiction.
\end{proof}

%DD
\begin{lemma} \label{L1}  Let $A$ be a maximal subdiagonal algebra.
If $h \in H^1$ is outer then
as an unbounded operator $h$ has dense range and trivial kernel.
Thus  $h = u |h|$ for a unitary $u \in M$.   Also, $\Phi(h)$ has
dense range and trivial kernel.
\end{lemma}

\begin{proof}  If $h$ is considered as an unbounded operator,
and if $p$ is the range projection of $h$, then since there exists a
sequence $(a_n)$ in $A$ with $h a_n \to 1$ in $L^1$-norm, we have
that $p^\perp = 0$. Thus the partial isometry $u$ in the polar
decomposition of $h$ is surjective, and hence $u$ is a unitary in $M$.
It follows that $|h|$ has dense range, and hence it, and $h$ also,
have trivial kernel.

For the last part note that $$L^1({\mathcal D}) = \Phi(H^1) =
\Phi([h A]_1) = [\Phi(h) {\mathcal D}]_1.$$
Thus we can apply the above arguments to $\Phi(h)$ too.
 \end{proof}

%dd
There is a natural equivalence relation on outers.  The following is
proved similarly to the classical case:

\begin{proposition} \label{L2}  If $h \in H^p$ is outer,
and if $u$ is a unitary in ${\mathcal D}$, then $h' = u h$ is outer
in $H^p$ too.  If $h, k \in H^p$ are outer,  then $|h| = |k|$ iff
there is unitary $u \in {\mathcal D}$ with $h = u k$.  Such $u$ is
unique.
\end{proposition}

As in the classical case, if $h \in H^2$ is outer, then $h^2$ is
outer in $H^1$.  Indeed one may follow the proof on \cite[p.\
229]{SW}, and the same proof shows that a product of any two outers
 is outer.
%In the non-antisymmetric case w
We do not know at the  time of writing whether every outer in $H^1$
is the square of an outer in $H^2$.

%dd
We now move towards a generalization of a beautiful classical
characterization of outer functions.  We will need Marsalli and
West's important earlier factorization result
generalizing the classical Riesz factorization. This states that 
%ll 
given any $1 \leq p, q, r \leq \infty$ with $1/p + 1/q = 1/r$, and 
$x \in H^r$, then $x = y z$ for some $y \in H^p$ and $z \in H^q$ 
\cite{MW}. (Marsalli and West actually only consider the case 
$r = 1$, but their proof generalises readily. A slightly sharpened 
version of their result may be found in \cite{LL1}.)

\begin{theorem}  \label{chr0}
If $A$ is an antisymmetric subdiagonal algebra  then
%dd
$h \in H^p$ is outer iff
$\Delta(h) = |\tau(h)| > 0$.
\end{theorem}

\begin{proof}
%dd
The case that $p = 2$ follows  exactly as in the proof of
\cite[Theorem 5.6]{Ho}. The case for general $p$ follows from the $p
= 1$ case and Lemma \ref{L3}. Hence we may suppose that $p = 1$. By
the Riesz factorization mentioned above,  $h = h_1h_2$ for $h_1, h_2
\in H^2$. Since $\Delta(h_1)\Delta(h_2) = \Delta(h) =
\Delta(\Phi(h)) = \Delta(\Phi(h_1))\Delta(\Phi(h_2)) > 0$ and
$\Delta(h_i) \geq \Delta(\Phi(h_i))$ for each $i = 1, 2$ (by the
generalized Jensen inequality), we must have $\Delta(h_i) =
\Delta(\Phi(h_i)) > 0$ for each $i = 1, 2$. Thus both $h_1$ and
$h_2$ must be outer elements of $[A]_2$, by the first line of the
proof. However we said above that a product of outers is outer.

Finally, if $h$ is outer then $1 = \lim_n \, h a_n$ where $a_n \in
A$, so that $1 = \lim_n \, \tau(h) \tau(a_n)$, forcing $\tau(h) \neq
0$. By the generalized Jensen inequality,  $\Delta(h) \geq
|\tau(h)|$. By Corollary \ref{vnc3}, for any $a_0 \in A_0$ we have
$$\Delta(h)
 \leq \tau(||h| (1 - a_0)|) =  \tau(|h - h a_0|). $$   Since
$h$ is outer, it is easy to see that $[h A_0]_1 = [A_0]_1$, and thus
we may replace $h a_0$ in the last inequality by $h - \tau(h)1 \in
[A_0]_1$. Thus $\Delta(h) \leq \tau(|\tau(h) 1|) = |\tau(h)|$.
\end{proof}
%LLe

{\bf Remarks.}  1) \  In the general non-antisymmetric case, one
direction of the last assertion of the last theorem is not quite
true as written. Indeed in the case that $A = M = L^\infty[0,1]$,
then outer functions in $L^2$ are exactly the ones which are a.e.\
nonzero. One can easily find an increasing function $h : [0,1] \to
(0,1]$,
% a step function with
%the step-lengths decreasing to $0$ as one approaches $0$,
which
%is outer  (so $1 \in [h L^\infty[0,1]]_2$) but
satisfies $\Delta(h) = 0$ or equivalently $\int_0^1 \, \log h = -
\infty$.
%$g_n(x) = e^n$ for $x \in [0,\frac{1}{n}]$, then clearly $g_n h \to
%1$
In the example of the upper triangular matrices ${\mathcal T}_2$ in
$M_2$, it is easy to find outers $h$ (i.e.\ invertibles in
${\mathcal T}_2$) with $0 < |\tau(h)| < \Delta(h)$.

%{\bf Remark.}
2) \ As in \cite[2.5.1 and 2.5.2]{SW},
%The proofs above contains
 $h$ is  outer in $H^2$  iff $|h| - \tau(h) u^* \in [|h|A_0]_2$, where
  $u$ is the
partial isometry in the polar decomposition of $h$, and iff $|h| -
\tau(h) u^*$ is the projection of $|h|$ onto $[|h|A_0]_2$. There may
be a gap in the proof of \cite[2.5.3]{SW}. The following principle
which they seem to be using is not correct:
%Following the classical method, we will use the following:
%\begin{lemma} \label{hilb}
If $K$ is a
subspace of a Hilbert space $H$, and if $\xi \notin K$, and $\xi -
\xi_0 \in H \ominus K$ with $\Vert \xi - \xi_0 \Vert = d(\xi,K)$,
then $\xi_0$ is the projection of $\xi$ onto $K$.
%\end{lemma}

\medskip

 As a byproduct of the method in the proof of the
 last theorem one
also has:

%ll
\begin{corollary}  Suppose that $A$ is antisymmetric.
An element in $H^r$ is outer iff it is the product of an outer in
$H^p$ and an outer in $H^q$, whenever  $1 \leq p, q, r \leq \infty$ 
with $1/p + 1/q = 1/r$.
\end{corollary}

\begin{corollary} \label{opp}  If
$A$ is antisymmetric, then $h$ is outer in
%dd
$H^p$ iff $[Ah]_p = H^p$.
 \end{corollary}

\begin{proof}
Replacing $A$ by $A^*$, we see that $\Delta(h) = \Delta(h^*) > 0$ is
equivalent to $h^*$ being outer in
%dd
$H^p(A^*) = (H^p)^*$. The latter
is equivalent to $(H^p)^* = [h^* A^*]_p$.  Taking adjoints again
gives the result.
\end{proof}

{\bf Remark.} \
%dd The same result for every $p \geq 1$ will follow  in
%exactly the same way from our main theorem below. Such
The last result has the consequence that the theory has a left-right
symmetry; for example our inner-outer factorizations $f = uh$ below
may instead be done with $f = h u$  (a different $u, h$ of course).
% (see e.g. \ the last part of the proof of
%Theorem \ref{chr}).

\medskip

Another classical theorem of Riesz-Szeg\"o states that if $f \in
L^1$ with $f \geq 0$, then $\int \, \log f \, > - \infty$ iff $f =
|h|$ for an outer $h \in H^1$ iff $f = |h|^2$ for an outer $h \in
H^2$. One may easily generalize such classical results if the
algebra $A$ is antisymmetric.

\begin{theorem}  \label{abssz}  If $f \in L^1(M)_+$ and
if $A$ is antisymmetric
then the following are equivalent:
\begin{itemize}
\item [(i)]  $\Delta(f) > 0$,
\item [(ii)]  $f \notin [f A_0]_1$,
\item [(iii)]  $f = |h|$ for $h$ outer in $H^1$.
\end{itemize}
\end{theorem}

\begin{proof}  (i) $\Leftrightarrow$ (ii) \ See Corollary \ref{vnc3}.

(ii) $\Rightarrow$ (iii) \ By the BN-factorization
\ref{BNpa} (b),  $f = u h$ for a unitary $u$ and an outer $h \in H^1$.
Since $f \geq 0$ we have $f = (f^* f)^{\frac{1}{2}} = |h|$.

(iii)  $\Rightarrow$ (ii) \ Conversely, if $h$ is outer then by the
%last part of
Lemma \ref{L3} we have $h \notin [h A_0]_1$.  By
%the other part of
Lemma \ref{L1} it follows that $|h| \notin [|h|
A_0]_1$.
\end{proof}

 %It is essentially a long outstanding open question as to
%whether (i) implies (iii) for general maximal subdiagonal algebras
%(see \cite[Remark after Theorem 8.1]{PX}).

We now turn to the topic of `inner-outer' factorizations. It can be
shown as in the classical case that if $f = uh$ is a
%`Beurling-Nevanlinna
factorization' of an $f \in L^p(M)$, for a unitary $u \in M$ and $h$
outer in $H^1$, then this factorization is unique up to a unitary in
${\mathcal D}$. Namely, if $u_1 h_1 = u_2 h_2$ were two such
factorizations, then $u_2 = u_1 u$ and $h_2 = u^* h_1$. See
\cite{BL6} for the easy details.

\begin{theorem}  \label{mas} If $f \in L^1(M)$, if $1 \leq p < \infty$,
 and if $A$ is antisymmetric,
then the following are equivalent:
\begin{itemize}
\item [(i)]  $\Delta(f) > 0$, \item [(ii)]  $f \notin [f A_0]_1$,
 \item [(iii)]  $|f|^{\frac{1}{p}} \notin [|f|^{\frac{1}{p}} A_0]_p$,
 \item [(iv)]  $f = uh$ for a unitary $u \in M$ and $h$ outer in $H^1$.
\end{itemize}
The factorization in {\rm (iv)} is unique up to a unimodular
constant.
\end{theorem}

\begin{proof}
(i) $\Rightarrow$ (iv) \ If $f \in L^1(M)$ with $\Delta(f) > 0$,
then $\Delta(|f|) > 0$, and so $|f| = |h|$ with $h$  outer in $H^1$,
by the previous result.  It follows by Lemma \ref{L1} that $f = u v
h$ for a partial isometry $u$ and a unitary $v$ in $M$. We have
$\Delta(f) = \Delta(u) \Delta(v) \Delta(h)$ by Theorem \ref{Fkd}
(4), so that $\Delta(u)
> 0$. Thus $\Delta(u^* u)
> 0$.  This forces $u^* u = 1$ by \cite[p.\ 606]{AIOA}, and so $u$
is unitary.

(iii) $\Leftrightarrow$ (i) \ See Corollary \ref{vnc3}.

(iv) $\Rightarrow$ (ii) \ If $f \in [fA_0]_1$, then $h \in [h
A_0]_1$, and by Lemma \ref{L3} we obtain a contradiction.

(ii) $\Rightarrow$ (iv) \  This is the BN-factorization \ref{BNpa} (b).

(iv) $\Rightarrow$ (i) \ We have $\Delta(u h) = \Delta(u) \Delta(h)
= \Delta(h)
%dd
> 0$, by previous results.

%DD
The uniqueness assertion follows from the remark above the theorem.
\end{proof}

{\bf Remarks.}  1) \ The $u$ in (iv) is necessarily in $[fA]_1$
(indeed if $h a_n \to 1$ with $a_n \in A$, then $f a_n = u h a_n \to
u$).

2) \ Suppose that in Theorem \ref{mas}, $f$ is also in $H^1$. Then,
first,  the $u$ in (iv) is necessarily in $[fA]_1 \subset H^1$, by
Remark 1. So $u \in H^1 \cap M = A$ (using $\tau$-maximality for
example).
 Thus $u$ is `inner' (i.e.\ is a unitary in $H^\infty = A$).  Second, note that (i)--(iv) will hold if $\tau(f)
\neq 0$. Indeed if $f \in H^1$ and if $f \in [f A_0]_1$ with $f a_n
\to f$ for $a_n \in A_0$, then $\tau(f) = \lim_n \, \tau(f a_n) =
0$.

\begin{corollary} \label{pth}  If
$A$ is antisymmetric and   $f \in L^1(M)_+$ then $\Delta(f) > 0$ iff
$f = |h|^p$ for an outer $h \in H^p$. \end{corollary}

\begin{proof}  ($\Rightarrow$) \  By
the previous result,  $f^{\frac{1}{p}} \notin [f^{\frac{1}{p}}
A_0]_p$, and
%the converse
%is fairly obviously true if $A$ is antisymmetric. Thus if $A$ is
%antisymmetric then,
so by the Beurling-Nevanlinna  factorization \ref{BNpa} (b) we have
$f^{\frac{1}{p}} = u h$, where $h$ is outer in $H^p$, and $u$ is
unitary.  Thus $f = (f^{\frac{1}{p}} f^{\frac{1}{p}})^{\frac{p}{2}}
=
 (h^* h)^{\frac{p}{2}}  = |h|^p$.
%, and $f^{\frac{1}{2}} = |h|$.

($\Leftarrow$) \ If $f = |h|^p$ for an outer $h \in H^p$ then
$\Delta(f) = \Delta(|h|)^p > 0$ by  Theorem \ref{Fkd} (3) and
%dd
Theorem  \ref{chr0}.
\end{proof}

\begin{corollary} \label{pth2}  If
$A$ is antisymmetric and   $f \in L^p(M)$ then $\Delta(f) > 0$ iff
$f = u h$ for a unitary $u$ and an outer $h \in H^p$.
%dd
If $p = 2$
and $v$ is the orthogonal projection of $f$ onto $[fA_0]_2$,
then $u$ is the partial isometry  in the polar decomposition
of $k-v$, and $h = u^* k$.
\end{corollary}

\begin{proof}  ($\Rightarrow$) \ By Theorem \ref{mas} we
obtain the factorization with outer $h \in H^1$. Since $|f| = |h|$
we have $h \in L^p(M) \cap H^1 = H^p$ (using \cite[Proposition 2]{Sai}).

($\Leftarrow$) \ As in the proof that (iv) implied (i) in
 Theorem \ref{mas}.

%dd
See \cite{BL6} for the proof of the last part (this is not used below).
\end{proof}

%dd
An obvious question is whether there are larger classes of
subalgebras of $M$ besides subdiagonal algebras for which such
classical factorization theorems hold. The following shows that,
with a qualification, the answer to this is in the negative. We omit
the proof, which may be found in \cite{BL6}, and proceeds by showing
that $A$ satisfies (ii) in Theorem \ref{B-L} above.

\begin{proposition} \label{gian}  Suppose that $A$ is a tracial
subalgebra of $M$, such that every $f \in L^2(M)$ with $\Delta(f) >
0$ is a product $f = u h$ for a unitary $u$ and an outer $h \in
[A]_2$.   Then $A$ is a finite maximal subdiagonal algebra.
\end{proposition}

%We will apply this lemma in the case that

{\em Question.}  \ Is there a characterization of outers in $H^1$ in
terms of extremals, as in the deLeeuw-Rudin theorem of e.g.\
\cite[p.\ 139--142]{Hobk}, or \cite[pp.\ 137-139]{Gam}?

\section{Logmodularity, operator spaces,
 and the uniqueness of extensions}

In the material above, we have not used operator spaces. This is not
because they are not present, but rather because they are not
necessary.  However, the algebras above do have interesting operator
space properties, and this would seem to add to their importance.
For example:

\begin{theorem} \label{3lmis2} {\rm  \cite{BL1,BLM} } \  Suppose that $A$ is
a logmodular subalgebra of a unital
$C^*$-algebra $B$, with $1_B \in A$.
\begin{enumerate}
\item [{\rm (1)}] Any unital completely contractive (resp.\
completely isometric) homomorphism $\pi  \colon A \rightarrow
B(H)$ has a unique extension to a completely positive and
completely contractive (resp.\ and completely isometric) map from
$B$ into $B(H)$.
\item [{\rm (2)}]  Every $*$-representation of $B$ is a boundary
representation for $A$  in the sense of {\rm \cite{Arv}}.
\end{enumerate}
\end{theorem}

\begin{proposition}  The $C^*$-envelope (or `noncommutative
Shilov boundary') of a maximal subdiagonal subalgebra of $M$ is $M$
again.  If $M$ is also injective then the injective envelope of a
maximal subdiagonal subalgebra of $M$ is $M$ again.
\end{proposition}

\begin{proof}  The proof of the first
assertion may be found in \cite{BL1,BLM}.
The last statement follows from the earlier one, since $M = I(M) =
I(C^*_e(A)) = I(A)$.   The last equality is valid for any unital
operator space $A$ by the `rigidity' characterization of the
injective envelope \cite[Section 4.2]{BLM}, since $A \subset
C^*_e(A) \subset I(A)$.
\end{proof}

There is a partial converse to some of the above:

\begin{theorem} \label{ex}   {\rm  \cite{BL4} } \ Suppose that
$A$ is a subalgebra of a unital  $C^*$-algebra $B$ such that $1_B \in A$, and
suppose that $A$ has the property that for every Hilbert space
$H$, every completely contractive
unital homomorphism $\pi : A \to B(H)$ has a unique completely
contractive (or equiv. completely positive) extension $B \to B(H)$.
Then $B = C^*_e(A)$, the $C^*$-envelope of $A$.
\end{theorem}

{\em Question:}  \  If $A$ is a tracial subalgebra of $M$ such that
every completely contractive
unital homomorphism $\pi : A \to B(H)$ has a unique completely
contractive extension $B \to B(H)$,
then is $A$ maximal subdiagonal?

\smallskip

If this were true, it would be a noncommutative analogue of Lumer's result from
\cite{Lum}.

\medskip

%DD
{\bf Closing remarks/further open questions.} \ 1)
 A question worthy of consideration is the
extent to which the results surveyed above  may be extended to the
setting of type III von Neumann algebras. Although most results in
Arveson's paper \cite{AIOA} are stated for subdiagonal subalgebras
of von Neumann algebras with a faithful normal tracial state, he
does also consider subalgebras of more general von Neumann algebras.
It would be interesting if there was some way to extend some of our
results to this context. See e.g.\
%D
recent work of Xu \cite{Xu}.
 %the question ofmaximality for
%on  this  class of algebras.
However there are some major challenges
to overcome if the full cycle of ideas presented above is to extend
to the type III case. We
 will try to elaborate this point. In the preceding theory the Fuglede-Kadison
 determinant plays a vital role. One expects that for a type III theory to work,
 a comparable quantity would have to found in that context. The problem is that
any von Neumann algebra which admits of a Fuglede-Kadison
determinant \emph{necessarily} admits of a faithful normal tracial
state (see e.g.\ \cite[Theorem 3.2]{Kai}).   Hence if the theory is to
extend to the case of, say, a type III algebra $M$ equipped with a
faithful normal state $\varphi$, then in order to have the necessary
tools at hand, one would first have to establish a theory of a
determinant-like quantity defined in terms of $\varphi$. Clearly
such a quantity cannot be a proper Fuglede-Kadison determinant as
such, but if it exists, one expects it to exhibit determinant-like
behavior with respect to the canonical modular
 automorphism group induced by $\varphi$: a kind of modular-determinant.
In support of the contention that such a quantity exists we note that
at least locally $\varphi$ does induce a determinant on $M$. That
is, any maximal abelian subalgebra of $M$ does admit of a
Fuglede-Kadison determinant, since the restriction of $\varphi$ to such a
 subalgebra is trivially a faithful normal tracial state. The challenge is
 to find a quantity which will yield a global expression for this local behavior.

%Telephone: +27-12-998-7447 (H) (Cell:)   +27-732865769
% Office ++27-12-429-6707

\medskip

2) A nice question communicated to us by Gilles Godefroy at this
conference is whether subdiagonal algebras have unique preduals.
This question is very natural in the light of his positive result in
this direction \cite[Theoreme 33]{Gode} in the function algebra case
(generalizing Ando's classical result on the uniqueness of predual
for $H^\infty(\Ddb)$).

\medskip

3) \ Another interesting question suggested to us at this conference, by
R. Rochberg, is whether the  generalization of the Helson-Szeg\"o
theorem to weak* Dirichlet algebras \cite{HiR}, and its corollary on
invertibility of Toeplitz operators, has a noncommutative variant in
this setting. This would probably need to use the real variable
%ll
methods developed by Marsalli and West in \cite{MW}, as well as their 
theory of Toeplitz operators with noncommuting symbols developed in 
\cite{MW2}.

\medskip

4) \  B. Wick has suggested to us that it might be worth
investigating noncommutative analogues of
\cite[Theorem 18.18]{Rud}, which may be viewed as 
a simple `uniform variant' of the corona theorem.

\medskip

5) \ Can one characterize (complete) isometries between
noncommutative $H^p$ spaces?  We have been informed by
Fleming and Jamison that there have been  
recent breakthroughs in the study 
of isometries between classical $H^p$ spaces. 
See \cite{LL1} for some work in this direction
in the noncommutative case.

\medskip

%B 
6) \  It seems extremely worthwhile,
and Arveson has also suggested
to us, to investigate certain subdiagonal algebras 
(or algebras which are not far from being subdiagonal)  
coming from free group examples as in \cite[Section 3]{AIOA}
in the framework of current free probability theory.
For example, it seems that aspects of our subject are not very far 
from some perspectives from the recent studies, by Haagerup and 
his collaborators, of the invariant subspace problem relative 
to a finite von Neumann algebra (see e.g.\ \cite{HS,HS2}).

\medskip

%BB

7) \ A most interesting project would be the widening 
of the class of algebras to which (parts of) the theory above 
may be extended.   There may in fact be several directions
in which to proceed.  For example, as in development 
of the commutative theory (see e.g.\ \cite{Gam,BK}), one could try to replace the 
requirement that $\Phi$ is multiplicative on $A$, by 
the multiplicativity of another conditional expectation $\Psi$.
If $\Delta(\Psi(a)) \leq \Delta(a)$ for $a \in A$, then one could view  
$\Psi$ as a {\em noncommutative Jensen measure}, and develop
a theory of the latter objects parallelling the important classical
theory.   Other clues to the enlarging of the class
of algebras might come from item 6 above, or from 
a closer examination of 
known classes of noncommutative
algebras which have some intersection with the 
class of subdiagonal algebras, and which do  
have noncommutative Hardy space properties.  
Such algebras have been studied in the 1970s and later,
beginning with \cite{Zs}.

\medskip

{\bf Acknowledgements.}
%dd
 We thank  Bill Arveson
for continual encouragement, and for several historical insights. We
thank  Barry Simon for relating to us some interesting mathematical
history of generalizations of Szeg\"o's theorem, and for pointing
out Verblunsky's precedence to the result usually attributed to
Kolmogorov and Krein.  
%B
We thank Gelu Popescu for interesting discussions. 
%B 
Finally, we are greatly indebted to Quanhua Xu for
extended conversations prompted by \cite{BL6}, and 
for many insightful and valuable comments.  He has recently
continued the work we did in \cite{BL6} by extending this
$H^p$ theory to values $0 < p < 1.$    
%BB
This, together with other very interesting related results
of his, should be forthcoming soon.


\begin{thebibliography}{99}
%B
\bibitem{PTAG}   W. B. Arveson, {\em Prediction theory and 
group representations,}  Ph. D. Thesis, UCLA, 1964.

 \bibitem{AIOA}  W. B. Arveson, Analyticity
in operator algebras, {\em Amer.\ J.\ Math.\ } {\bf 89} (1967),
578--642.

\bibitem{Arv}  W. B. Arveson,
 Subalgebras of $C^{*}$-algebras,
{\em  Acta Math.} {\bf 123 }(1969), 141-22.

%BB
\bibitem{BK}  K. Barbey and H. K\"onig, 
{\em Abstract analytic function theory and Hardy algebras,}
 Lecture Notes in Mathematics, Vol.\ 593,  
Springer-Verlag, Berlin-New York, 1977.

\bibitem{BL1}    D. P. Blecher and L. E. Labuschagne,
Logmodularity and isometries of operator algebras, {\em
 Trans.\ Amer.\ Math.\ Soc.\ } {\bf 355}  (2003),  1621--1646.

\bibitem{BL2}    D. P. Blecher and L. E. Labuschagne,
Characterizations of noncommutative $H^\infty$,   {\em
 Integr.\ Equ.\ Oper.\ Theory}  {\bf 56} (2006), 301-321.
% : Online First (2006), DOI 10.1007/s00020-006-1425-5.
 %dd {\bf 355}  (2006),  1621--1646.

\bibitem{BL3}    D. P. Blecher and L. E. Labuschagne,
A Beurling theorem for noncommutative $L^p$, To appear, {\em J.\
Operator Theory}.

\bibitem{BL4}    D. P. Blecher and L. E. Labuschagne,
Noncommutative function theory and unique extensions, To appear,
{\em Studia Math.}  

\bibitem{BL6}    D. P. Blecher and L. E. Labuschagne,
Applications
of the Fuglede-Kadison determinant: 
Szeg\"o's theorem and outers for noncommutative $H^p$, 
(ArXiv:  math.OA/0609662).  

\bibitem{BLM}  D. P. Blecher
and C.  Le Merdy, {\em Operator algebras and their modules---an operator space
approach,} Oxford Univ. Press, Oxford, 2004.

\bibitem{Brme}  L. G. Brown, {\em Lidski\u\i's theorem in the
type ${\rm II}$ case,}
Geometric methods in operator algebras (Kyoto, 1983),  1-35, Pitman
Res. Notes Math. Ser., 123, Longman Sci. Tech., Harlow, 1986.



\bibitem{E}  R. Exel, Maximal subdiagonal algebras,
{\em Amer.\ J.\ Math.\ } {\bf 110} (1988), 775-782.


\bibitem{E2}  R. Exel,  The F. and M. Riesz theorem for $C^*$-algebras,
   {\em J. Operator Theory} {\bf   23}  (1990),   351-368.



\bibitem{FK}  T. Fack and H. Kosaki, Generalized s-numbers of $\tau$-measurable
operators, {\em Pacific J. Math} {\bf123} (1986), 269-300.

\bibitem{Gam}   T. W. Gamelin,
 {\em Uniform Algebras,} Second edition, Chelsea, New
York, 1984.

\bibitem{Gode} G. Godefroy,   Sous-espaces bien disposes de
$L\sp{1}$-applications. {\em Trans.\ Amer.\ Math.\ Soc.} {\bf  286}
  (1984),   227-249.

\bibitem{HS}   U. Haagerup and H. Schultz,  Brown measures of unbounded operators affiliated
with a finite von Neumann algebra,   Preprint (2006)
math.OA/0605251

%B
\bibitem{HS2}   U. Haagerup and H. Schultz,  
Invariant subspaces for operators in a general $II_1$-factor, Preprint (2006)
math.OA/0611256 

%dd
\bibitem{Hel}  H. Helson, {\em Lectures on invariant subspaces,}
 Academic Press, New York-London, 1964. 

\bibitem{HL}  H. Helson 
and D. Lowdenslager, Prediction theory and Fourier series
in several variables, {\em Acta Math.} {\bf 99} (1958), 165-202.

 \bibitem{HL2}  H. Helson and D. Lowdenslager, 
Prediction theory and Fourier 
series in several variables. II,  {\em Acta Math.} {\bf 106} (1961), 175--213.

\bibitem{HiR}  I. I. Hirschman, Jr. and R. Rochberg, Conjugate function theory in
weak* Dirichlet algebras, {\em J.\ Funct.\ Anal} {\bf 16}
(1974), 359-371.

\bibitem{Hobk}  K. Hoffman, {\em Banach spaces of analytic functions},
 Dover (1988).



\bibitem{Ho}  K. Hoffman, Analytic functions and logmodular
Banach algebras, {\em  Acta Math.} {\bf 108} (1962), 271-317.

\bibitem{HR}  K. Hoffman and H.  Rossi, Function theory and
multiplicative linear functionals, {\em Trans.\ Amer.\ Math.\
Soc.} {\bf 116} (1965), 536-543.

\bibitem{JOS} G.  Ji, T.  Ohwada, K-S. Saito,
    Triangular forms of subdiagonal algebras,
{\em Hokkaido Math.\ J.} {\bf 27}, (1998), 545-552.



\bibitem{JS}  M. Junge and D. Sherman, Noncommutative $L^p$-modules,
{\em J.\ Operator Theory} {\bf 53} (2005), 3-34.



\bibitem{KR}   R. V. Kadison and J. R. Ringrose,
{\em Fundamentals of the theory of operator algebras,} Vols. 1 and 2,
 Graduate Studies in Mathematics,
Amer.\ Math.\ Soc.\ Providence, RI, 1997.

\bibitem{Kai}  S. Kaijser,
 On Banach modules II. Pseudodeterminants and traces, {\em
  Math. Proc. Cambridge Philos. Soc.} \ {\bf 121}  (1997),  325--341.

%B
\bibitem{KT} S.  Kawamura and J.  Tomiyama,
 On subdiagonal algebras associated with flows in operator algebras,
{\em   J. Math. Soc. Japan} {\bf  29}  (1977),  73--90.

\bibitem{LL1}  L. E. Labuschagne,  Analogues
of composition operators on non-commutative $H^p$-spaces,
{\em J. Operator Theory} {\bf 49}  (2003),   115-141.

\bibitem{LL2} L. E. Labuschagne,  Composition operators on
 non-commutative $L_p$-spaces, {\em Expo. Math.} {\bf 17} (1999),
 429-468.

\bibitem{LL3}  L. E. Labuschagne,  A noncommutative Szeg\"o theorem
for subdiagonal subalgebras of von Neumann algebras, {\em Proc.
Amer. Math. Soc.}, {\bf 133} (2005), 3643-3646.


%B
\bibitem{LM}  R. I. Loebl and P. S. Muhly,  Analyticity and flows
in von Neumann algebras, {\em   J.\ Funct.\ Anal.\ } {\bf 29} (1978),
214--252.


\bibitem{Lum} G. Lumer,  Analytic functions and Dirichlet
problems, {\em  Bull. Amer. Math. Soc.} {\bf 70} (1964), 98-104.

\bibitem{MMS}  M. McAsey, P. Muhly, and K.-S. Saito,  Nonselfadjoint
crossed products (invariant subspaces and maximality),
{\em  Trans. Amer. Math. Soc.} {\bf 248} (1979), 381-409.

\bibitem{MW}  M. Marsalli and G. West,  Noncommutative
$H^p$-spaces,  {\em J.\ Operator Theory} {\bf 40} (1998), 339--355.

%ll
\bibitem{MW2}  M. Marsalli and G. West,  Toeplitz operators with 
noncommuting symbols,   {\em  Integr.\ Equ.\ Oper.\ Theory} {\bf 32} 
(1998), 65 -- 74.

\bibitem{MS}  P. S. Muhly and B. Solel,
Hardy algebras, $W\sp *$-correspondences and interpolation theory,
{\em  Math.\ Ann.}  {\bf 330}  (2004),  353-415.

\bibitem{Nak}  T. Nakazi, Extended weak-$*$ Dirichlet algebras,
\textit{Pacific
%Invariant subspaces of weak-$*$ Dirichlet algebras. \textit{Pacific
J. Math.} {\bf 81}(1979), 493-513.


\bibitem{N}  T. Nakazi and Y.  Watatani, Invariant
subspace theorems for subdiagonal algebras, {\em J.\
Operator Theory} {\bf
37}  (1997),   379-395.

\bibitem{Neu}  B. H. Neumann, On ordered groups. {\em Amer. J. Math.}
{\bf  71}  (1949), 1--18.

\bibitem{Nik}  N.K. Nikolski
\textit{Operators,
functions, and systems: an easy reading. Vol. 1. Hardy, Hankel,
and Toeplitz.}
 Mathematical Surveys and Monographs, 92. American
 Mathematical Society, Providence, RI, 2002.


\bibitem{Nel} E. Nelson, Notes on noncommutative integration,
{\em J. Funct. Anal.} {\bf 15} (1974), 103-116.


\bibitem{PX}   G. Pisier and Q. Xu, {\em Noncommutative $L^p$ spaces,}
in   Vol 2. of
Handbook on Banach spaces, Ed. W. B. Johnson and J. Lindenstrauss,
North-Holland, Amsterdam, 2003.

\bibitem{Pitts}  D. Pitts, Factorization problems for nests:
factorization methods and characterizations of the universal
factorization property,
{\em J. Funct. Anal.} {\bf 79} (1988), 57-90.

\bibitem{Pop}  G. Popescu,
Entropy and multivariable interpolation,  Mem.\
 Amer.\ Math.\ Soc.  {\bf 184}(2006), No.\ 868.



\bibitem{Po}  S. C. Power,   Factorization in analytic operator
algebras, {\em J. Funct. Anal.} {\bf 67} (1986), 413-432.

\bibitem{Rand}  N.
Randrianantoanina,
Hilbert transform associated with finite maximal subdiagonal algebras,
{\em   J.\ Austral.\ Math.\ Soc.\ Ser.\ A}  {\bf 65}  (1998),  388-404.

\bibitem{Rud}  W. Rudin, {\em Real and complex analysis,} Third Ed.,
 McGraw-Hill (1987).

 \bibitem{Sai}   K.-S. Saito,  A note on invariant subspaces for
finite maximal subdiagonal algebras,  {\em  Proc.\ Amer. Math. Soc.} {\bf 77}
 (1979), 348-352.

\bibitem{Sak} S. Sakai, {\em $C^*$-algebras and $W^*$-algebras}, Springer,
New York, 1971.

\bibitem{SW}  T. P.  Srinivasan and J-K. Wang, {\em Weak*-Dirichlet algebras,}
 In {\em Function algebras}, Ed. Frank T. Birtel,
 Scott Foresman and Co., 1966, 216-249.

\bibitem{Tak}   M. Takesaki, {\em Theory of operator algebras I},
Springer, New York, 1979.

\bibitem{Tak2}   M. Takesaki, {\em Theory of operator algebras II},
Springer, New York, 2003.

\bibitem{Terp}  M. Terp, {\em $L^p$ spaces associated with von Neumann algebras,}
 Notes, Math.\ Institute, Copenhagen Univ.  1981.


\bibitem{Ver} S. Verblunsky, On positive
harmonic functions (second paper), {\em
Proc.\ London Math.\ Soc.} \textbf{40} (1936), 290-320.

\bibitem{Xu}  Q. Xu,  On the maximality of subdiagonal algebras,
 {\em J. Operator Theory}   {\bf 54}  (2005), 137-146.

\bibitem{Zs}  L. Zsid\'o,
Spectral and ergodic properties of the analytic generators,
{\em J. Approximation Theory} {\bf 20} (1977), 77-138.
\end{thebibliography}
\end{document}